\documentclass[runningheads]{llncs}
\usepackage{graphicx}

\usepackage{amssymb,amstext,amsmath,mathabx}

\usepackage[utf8]{inputenc}
\usepackage{hyperref}
\usepackage{verbatim}
\usepackage{mathtools}
\usepackage{color}
\usepackage[normalem]{ulem} 
\usepackage[mathscr]{euscript}
\usepackage{wasysym}
\newcommand{\no}[1]{\widebar{#1}}

\renewcommand{\mid}{|}
\def\F{\mathcal F}
\def\P{\mathcal P}
\def\B{\mathcal B}

\def\I{\mathcal I}
\def\H{\mathcal H}
\def\pr{\mathbb{P}}
\def\prev{\mathbb{P}}
\def\L{\mathcal{L}}
\def\K{\mathcal{K}}

\def\G{\mathcal{G}}

\def\S{\mathcal{S}}
\def\B{\mathcal{B}}

\def\C{\mathscr{C}}
\def\D{\mathscr{D}}

\title{
Subjective probability, trivalent logics and compound conditionals 
}
\titlerunning{Subjective probability, trivalent logics and compound conditionals }
\author{Angelo Gilio$^1$ \and
	Giuseppe Sanfilippo$^2$}
\institute{Dept. SBAI, University of Rome ``La Sapienza'', Rome, Italy (retired)\\
	\email{angelo.gilio@sbai.uniroma1.it}\\
	\and
	Dept. Mathematics and Computer Science, University of Palermo, Palermo, Italy \\
	\email{giuseppe.sanfilippo@unipa.it}}

\begin{document}
	\maketitle

	\begin{abstract}
In  this work we first illustrate the subjective theory of de Finetti. We recall the notion of coherence for both the  betting scheme and the penalty criterion, by considering the  unconditional and conditional cases.
We show  the equivalence of the two criteria by   giving  the geometrical interpretation of coherence. We also consider the notion of coherence based on  proper scoring rules.  We discuss  conditional events in the trivalent logic of de Finetti and the  numerical representation of  truth-values. We check the  validity 
of selected  basic logical and probabilistic properties for  some  trivalent logics:
Kleene-Lukasiewicz-Heyting-de Finetti; Lukasiewicz; Bochvar-Kleene; Sobocinski.
We verify  that none  of these logics  satisfies  all the  properties. Then, we consider
  our approach to conjunction and disjunction of conditional events  in the setting of  conditional random quantities.  We verify that  all the basic logical and probabilistic properties (included the Fréchet-Hoeffding bounds)  are preserved in our approach. We also recall the characterization of  p-consistency and p-entailment by our notion of  conjunction.
	\end{abstract}

	\section{Introduction}
Conditionals play a relevant role in many fields,
such as philosophy, psychology of the uncertain reasoning, artificial intelligence, decision making, nonmonotonic reasoning, knowledge representation. In a pioneering paper de Finetti (\cite{deFi35}) introduced conditional events as tri-events and proposed a  three-valued logic.
Many  authors have  given important contributions to research on   conditionals;
see, e.g.,  
\cite{adams75,CoSV13,CoSV15,Douven19,DuPr94,edgington95,FlGH20,gilio02,GiSa14,GiSa19,hailperin96,Jeff91,Kauf09,Isberner2001,McGe89,Miln97,Mura11,Over21,PfSa17,SPOG18,StJe94,vanFraassen1976}.	
Usually, conjunction and disjunction of conditionals have been defined as suitable conditionals in the setting of trivalent logics; see, e.g., \cite{Baratgin18,benferhat97,Cala87,CiDu13,EgRS20,GoNW91,NgWa94,Scha68}. 
But, in this way,  basic logical and probabilistic properties valid for unconditional events are not preserved. For instance, the lower and upper bounds on the probability of the  conjunction, or the disjunction, of conditional events do not coincide with the  Fr\'echet-Hoeffding bounds. 
A more general approach to compound conditionals, 
 where the result of conjunction or disjunction is no longer a  three-valued object has been given in  \cite{Kauf09,McGe89}.     A related theory
 has been  developed in the setting of coherence (\cite{GiSa14,GiSa19,GiSa20,GiSa21}), by defining  conjoined and disjoined conditionals
  as suitable  {\em conditional random quantities}, with values in the interval $[0,1]$.   For some applications see, .e.g., \cite{GiPS20,SGOP20,SPOG18}. A further development of this approach to general compound conditionals has been given in \cite{FGGS22}. 
  Among various approaches to probability, we use the subjective interpretation of de Finetti which is  based on the coherence principle.  The theory of de Finetti  allows a flexible and gradual probabilistic treatment of uncertainty. 
Within this theory we can directly  assign conditional probabilities on arbitrary families of conditional events, without requiring algebraic structures, and  we can coherently extend   probabilities to further conditional events (\emph{fundamental theorem of probability} of de Finetti and its generalizations). A further advantage  is the possibility of properly managing conditioning events with zero probability. \\
We list below the main aspects considered in the paper.
\begin{itemize}
	\item[$\bullet$] We recall  the notion of coherence within the betting scheme, the penalty criterion, and proper scoring rules.
	\item[$\bullet$] We recall the geometrical interpretation of  coherence.
	\item[$\bullet$] We discuss in the betting framework the   numerical representation of a  conditional event $A|H$, by showing that the natural counterpart of the logical values  \emph{true}, \emph{false}, and \emph{void} are  the numerical values $1$, $0$, $P(A|H)$ of a suitable conditional random quantity (the indicator of $A|H$).   
    \item[$\bullet$] We consider some basic logical and probabilistic properties  valid in the case of unconditional events. Then we examine four trivalent logics and we verify that none of them satisfy all the properties. 
    \item[$\bullet$]  We illustrate the approach to compound conditionals,  where conjunctions and disjunctions are defined as suitable conditional random quantities.  We verify that in our approach all the basic (and other further) properties are preserved.    
\end{itemize}
The paper is organized  as follows. 	In Section \ref{SEC:COHERENCE}  we recall the notion of coherence of probability assessments of unconditional events, for both 
the  betting scheme and the penalty criterion. Then,  we recall the geometrical characterization of coherence given by de Finetti.  In \ref{SEC:CONDCOHERENCE}  we discuss conditional events in the trivalent logic of de Finetti and  the numerical representation of the truth-values, by showing that the numerical counterpart of the logical value \emph{void} for a conditional event $A|H$ is given by  $P(A|H)$.
Then we illustrate the  notion of coherence of conditional probability assessments for 
the  betting scheme and the penalty criterion. We verify  the equivalence of the two criteria by showing  that the   geometrical interpretation of coherence is the same in both criteria. 
 We  also briefly discuss the notion of coherence for conditional probability assessments  by means of a generic proper scoring rules. 
In Section \ref{SEC:TRIVALENT} we discuss compound conditionals in the setting of trivalent logics.  We consider four notion of conjunctions:
Kleene-Lukasiewicz-Heyting-deFinetti conjunction ($\wedge_K$), Lukasiewicz conjunction  ($\wedge_L$);  Bochvar internal conjunction, or Kleene weak conjunction ($\wedge_B$);  Sobocinski conjunction, or quasi conjunction ($\wedge_S$). By De Morgan Law we also consider the associated disjunctions $\vee_K,\vee_L,\vee_B,\vee_S$. 
Then, we examine some basic logical and probabilistic properties, valid for unconditional events,  and we verify that none of the previous trivalent logics satisfies all the properties.
In Section \ref{SEC:CONDRAND}  we recall our approach to the  conjunction and the disjunction of conditional events,  defined in the setting of coherence as suitable conditional random quantities. We verify that, when considering conditional events,  all the basic logical and probabilistic properties valid for unconditional events are preserved.  Finally, in Section \ref{SEC:FURTHER} we consider some further properties which are satisfied in our approach to compound  conditionals.
	\section{Coherent probability assessments}
\label{SEC:COHERENCE}
	Uncertainty is  present in almost all real problems, where we have to face  facts which will turn out to be true or false.  We describe uncertain facts, termed events, by non ambiguous logical  propositions.  
	Then, an event   $E$ is a two-valued logical entity which can be \emph{true}, or \emph{false}. The indicator of $E$, denoted by the same symbol, is  1, or 0, according to  whether $E$ is true, or false.  Thus, a symbol like $sE$ represents the product of the quantity $s$ and the indicator of the event $E$.
	The sure event and impossible event are denoted by $\Omega$ and  $\emptyset$, respectively. 
	Given two events $E_1$ and $E_2$,  we denote by $E_1\land E_2$, or simply by $E_1E_2$, (resp., $E_1 \vee E_2$) the logical conjunction (resp., the  logical disjunction).   
	The  negation of $E$ is denoted $\no{E}$.  We simply write $E_1 \subseteq E_2$ to denote that $E_1$ logically implies $E_2$, that is  $E_1\no{E}_2=\emptyset$.
	We recall that  $n$ events $E_1,\ldots,E_n$ are said to be  {\emph logically independent} when the number $m$ of constituents, or possible worlds, generated by them  is $2^n$, that is when the all the $2^n$ conjunctions $E_1^*\cdots E_n^*$, with $E_i^*\in\{E_i,E_i^*\}$, $i=1,2,\ldots, n$, are not impossible.

	In the subjectivistic  theory of de Finetti, the probability $P(E)$ that an individual  attributes to  an event $E$,  in his current state of knowledge, is a measure of  its degree of belief on  $E$ being true. In order to operatively assess  probabilities, de Finetti proposed  two (equivalent) criteria: $(i)$ \emph{betting scheme}; $(ii)$ \emph{penalty criterion}. 
	In both criteria, the consistency of the probability assessments is  guaranteed by a  coherence principle. All probabilistic properties follow  from coherence; moreover, if a probability  assessment does not satisfy some   basic probabilistic property, then it is not coherent.
	\subsection{Betting Scheme}
Within  this criterion, to assess $P(E)=p$ means that, for every real number $s$, you are willing to pay an amount $sp$ and to receive $sE$, that is to receive $s$, or $0$, according to whether $E$ is true, or $E$ is false, respectively. The random gain is $G = s(E-p)\in\{s(1-p),-ps\}$. The assessment $p$ is said to be \emph{coherent} if, for every $s$, it does not happen that the  values of the random gain are both positive, or both negative. Coherence requires that $P(E)\in[0,1]$, with $P(\Omega)=1$ and $P(\emptyset)=0$. More in general, given  an arbitrary family of events $\K$, a  probability assessment  on the events in  $\K$ is made by specifying a real function $P : \; \mathcal{K}
	\, \rightarrow \, \mathcal{R}$.
	For each finite  subfamily 
	$\mathcal{F} = \{E_1, \ldots, E_n\}$ of $\mathcal{K}$
	we denote by  $\mathcal{P} =(p_1, \ldots, p_n)$, where  $p_i = P(E_i) \, ,\;\; i \in \{1,\ldots,n\}$, 
	the  probability assessment on $\F$.  
	Within the betting scheme,  for  any arbitrary  real numbers $s_1, \ldots, s_n$,  
	you agree to pay $\sum_{i=1}^n s_ip_i$, by receiving the random amount $\sum_{i=1}^n  s_iE_i$.
	The random gain $G$, associated with the 
	pair $(\mathcal{F}, \mathcal{P}$) is the difference  between the  (random) amount that you receive and the amount that you pay, that is 
	$G =\sum_{i=1}^ns_iE_i-\sum_{i=1}^ns_ip_i =\sum_{i=1}^n s_i(E_i - p_i)$.
	The random quantity $G$ represents the net gain from engaging each transaction $(E_i - p_i)$, the scaling and meaning (buy or sell) of the transaction being specified by the magnitude and the sign of $s_i$, respectively. We observe that, by expanding the second member of the equality $\Omega=(E_1 \vee \no{E}_1) \wedge \cdots \wedge (E_n \vee \no{E}_n)$, we obtain the disjunction of  $2^n$ conjunctions $A_1 \cdots A_n$, where $A_i \in \{E_i, \no{E}_i\},\, i=1,\ldots,n$. By discarding the conjunctions which coincide with the impossible event $\emptyset$, the sure event $\Omega$ can be represented as a suitable disjunction $\Omega = C_1 \vee  \cdots \vee C_m$, where  $C_1, \ldots, C_m$ are the constituents, or possible cases, generated by $\F$, with $m \leq 2^n$. In the case of logical independence  $m=2^n$.
	\begin{example}
		The constituents associated with  $\F=\{A,B\}$, where  $A\subseteq B$,  are $C_1=AB$, $C_2=\no{A}B$ and $C_3=\no{A}\;\no{B}$.
	\end{example}
	
	Denoting by $g_h$ the value of $G$ when $C_h$ is true, it holds that $G\in \{g_1,\ldots,g_m\}$, with $g_h=\sum_{i: C_h \subseteq E_i} s_i - \sum_{i=1}^ns_ip_i$. 
	Then, within the  {\em betting  scheme}, we
	have
	\begin{definition}\label{DEF:COER-BET}  The function $P$ defined on $\mathcal{K}$ is said to be {\em coherent}
		if and only if, for every integer $n$, for every finite subfamily $\mathcal{F}$ of 
		$\mathcal{K}$ and for every  real numbers $s_1, \ldots, s_n$, one has:
		$\min  G \leq 0 \leq \max G$.
	\end{definition}
	The condition $\min  G \leq 0 \leq \max G$ can be equivalently written as: $\min  G \leq 0$, or  $\max  G \geq 0$.  As shown by Definition \ref{DEF:COER-BET}, a probability assessment is coherent if and only if, in any finite combination of $n$ bets, it does not happen that the values $g_1, \ldots, g_m$ are all positive, or all negative ({\em no Dutch Book}). 
	\begin{remark}\label{REM:GAIN}
		Given a probability assessment $\P=(p_1,\ldots,p_n)$  on a finite family $\F=\{E_1,\ldots,E_n\}$,  the coherence of $\P$ with Definition \ref{DEF:COER-BET}  means that, for every  real numbers $s_1, \ldots, s_n$,  the condition    $\min  G \leq 0 \leq \max G$  is satisfied. 
	\end{remark}	
	
	\begin{example}\label{Ex:coerenza1}
		Let $\P=(p_1,p_2,p_3)=(0.4,0.3,0.8)$ be a probability assessment on  $\F=\{E_1,E_2,E_3\}$, where $E_1=A$, $E_2=B$ and $E_3=A\vee B$, with $A,B$ logically independent.
		We show that the assessment $\P$ on $\F$ is incoherent.
		We set $s_1=s_2=-s_3=1$. The associated random gain is
		\[
		G=(A-0.4)+(B-0.3)-(A\vee B-0.8).
		\]
		The constituents associated with $\F$ are $C_1=AB$, $C_2=A\no{B}$, $C_3=\no{A}B$ and $C_4=\no{A}\;\no{B}$.
		The possible values  $g_1,\ldots,g_4$ of the random gain $G$, associated with $C_1,\ldots, C_4$, are
		\[
		\begin{array}{cl}
			g_{1}= & (1-0.4)+(1-0.3)-(1-0.8)=.0.6+0.7-0.2=1.1,\\
			g_{2}= & (1-0.4)+(-0.3)-(1-0.2)=0.6-0.3-0.2=0.1,\\
			g_{3}= & (-0.4)+(1-0.3)-(1-0.2)=-0.4+0.7-0.2=0.1,\\
			g_{4}= & (0-0.4)+(0-0.3)+(1-0.2)=-0.4-0.3+0.8=0.1. 
		\end{array}
		\]		
		As the possible value of $G$ are all positive, it follows that $\P$ is incoherent, indeed in this case $P(A\vee B)>P(A)+P(B)$, while coherence requires that $\max\{P(A),P(B)\}\leq P(A\vee B)\leq \min\{P(A)+P(B),1\}$.
	\end{example}

	\subsection{Penalty Criterion} 
	The notion of coherence can be also defined by exploiting a penalty criterion, which is equivalent to the betting scheme. With the assessment $\mathcal{P} =(p_1, \ldots, p_n)$ on $\mathcal{F}=\{E_1,\ldots,E_n\}$ we associate a random loss $\L=\sum_{i=1}^n (E_i-p_i)^2$, which represents the square of the distance between the (binary) random point $(E_1,\ldots,E_n)\in\{0,1\}^n$ and the prevision point $\P$. Then, within the  {\em penalty criterion}, we
	have
	\begin{definition}\label{DEF:COER-PENALTY}  The function $P$ defined on $\mathcal{K}$ is said to be {\em coherent} if and only if, for every integer $n$, for every finite subfamily $\mathcal{F}=\{E_1,\ldots,E_n\}$ of  $\mathcal{K}$,
		denoting by $\P=(p_1,\ldots,p_n)$  the restriction of $P$ to $\F$, 
		there does not exist $\mathcal{P}^* =(p_1^*, \ldots, p_n^*)$ such that: 
		$\L^* < \L$, where $\L$ and $\L^*$ are the losses associated with $\P$ and $\P^*$, respectively.
	\end{definition} 
	{\em Geometrical interpretation of coherence.} In the probabilistic theory of de Finetti coherence is also characterized by a geometrical interpretation.  We represent each constituent $C_h$ by a point $Q_h=(q_{h1},\ldots,q_{hn})$, where (for each index $i$) $q_{hi}=1$, or $q_{hi}=0$, according to whether $C_h \subseteq E_i$, or $C_h \subseteq \no{E}_i$, respectively. Notice that $q_{hi}$ is nothing but the value assumed by the indicator of $E_i$ when $C_h$ is true. 
	We  denote by $\I$ the convex hull of the points $Q_1,\ldots,Q_m$.  
	\begin{remark}\label{REM:LOSS}
		Given a probability assessment $\P=(p_1,\ldots,p_n)$  on a finite family $\F=\{E_1,\ldots,E_n\}$,  the coherence of $\P$ with Definition \ref{DEF:COER-PENALTY}  means that there does not exist $\mathcal{P}^* =(p_1^*, \ldots, p_n^*)$ such that 
		$\mathcal{L}^* < \mathcal{L}$; coherence of $\P$ is also  equivalent to the geometrical condition
		$\P\in\I$ (\cite{deFi70}). 
	\end{remark}	
	We will show later  
	(see Remark \ref{REM:EQUIV-CONV-GUAD}) that 
	\begin{equation}\label{EQ:CONV-GUAD}
	\min  G \leq 0 \leq \max G \Longleftrightarrow \;\;\; \P \in \I \,.
	\end{equation}
	Thus, by Remarks \ref{REM:GAIN} and \ref{REM:LOSS}, in order to determine the set of coherent assessments on $\F$, Definitions \ref{DEF:COER-BET} and  \ref{DEF:COER-PENALTY} are equivalent.
	
	We recall that, if the function $P$ defined on $\K$ is coherent, then it satisfies all the basic probabilistic properties. Moreover, it can be extended to a  finitely additive probability $P'$ on every algebra $\mathcal{A}$, with $\K\subseteq \mathcal{A}$.

	\section{Coherent conditional probability assessments}
	\label{SEC:CONDCOHERENCE}
	A more general situation is represented by the case where we want to assess a conditional probability $P(A|H)$, which represents the degree of belief on the event $A$, by assuming true the event $H$ (and nothing more). $P(A|H)$ is also called the probability of "$A$ given $H$", or the probability of the conditional event $A|H$, where $H \neq \emptyset$. \\ 
	Notice that a large number of philosophers and psychologists, in their research on conditionals, assume valid the {\em Equation} (see, e.g., \cite{douven11b,edgington95,Stal70}), or {\em Conditional Probability Hypothesis (CPH)}   (\cite{Cruz20,Over21,SGOP20}), by judging that the probability of a conditional "\emph{if $H$ then $A$}" coincides with the conditional probability $P(A|H)$. 
	\subsection{Numerical representation of truth values}
	In the theory of de Finetti (\cite{deFi35}) the conditional event $A|H$ is looked at as a three-valued logical entity, with values {\em true}, or {\em false}, or {\em void}, according to whether $AH$ is true, or $\no{A}H$ is true, or $\no{H}$ is true. \\ {\em What is an appropriate way of defining the indicator of $A|H$, i.e. to represent numerically the truth-values of $A|H$?} \\
	We observe that, likewise the unconditional case, the numerical counterparts of \emph{true} and \emph{false} are $1$ and $0$, respectively.  Concerning the third value \emph{void}, in the coherence-based approach of de Finetti (\cite{deFi35}, see also \cite{definetti95}) the evaluation of $P(A|H)$ is operatively based on the following {\em conditional bet}: \\
	 (Scheme 1). You evaluate $P(A|H)$, say $P(A|H)=x$,  when $H$ is uncertain. After $H$ is verified, the bet has effect  and you accept to pay the amount $sx$ in order to bet on $A$, by receiving  $s$ if $A$ is true, or 0 if $A$ is false. In the case where $H$ is not verified the bet has no effect. Within this conditional bet, we can say that \textit{if $H$ then I bet that $A$}, that is a bet on $A$ conditionally on $H$ being true. \\
	Coherence of the assessment $P(A|H)=x$ is checked by (only) considering the cases where the bet is effective, that is  $H$ is verified. The random gain  when the bet has effect (i.e., when $H$ true) is  $G_H=s(A-x)\in\{s(1-x),-sx\}$. Then, $x$ is coherent if and only if $\min G_H \leq 0\leq \max G_H$, that is $x\in[0,1]$, with $x=1$, or $x=0$, when  $AH=H$, or $AH=\emptyset$, respectively. \\
	The previous conditional bet   can be expressed in the following \textit{equivalent form}: \\
	(Scheme 2). If you assess $P(A|H)=x$, then (before knowing the truth value of $H$) you accept to pay $sx$, by receiving $s$ if $AH$ is true, or 0 if $\no{A}H$ is true, or $sx$ if $\no{H}$ is true (the bet is called off).  In this scheme $G=sH(A-x)$, with $G \in\{s(1-x),-sx,0\}$ and $G_H\in \{s(1-x),-sx\}$.
	By definition, for   checking  coherence only the cases in which the bet is not called off are considered.  Thus, coherence of $x$ amounts to $\min G_H \leq 0\leq \max G_H$, that is   
	\[
	\min \, \{s(1-x),-sx\} \; \leq \; 0 \; \leq \; \max \, \{s(1-x),-sx\} \,.
	\]
	\subsection{On the third value of conditional events}\label{THIRD}
	Focusing on the numerical counterpart of the logical value \emph{void} we recall that, in a bet on an unconditional event $A$, you pay for instance $P(A)$ by receiving (the value of the indicator) $A$.  Likewise, in the conditional case  you pay for instance $P(A|H)$ by receiving $AH+P(A|H)\no{H}$, which should represent the indicator of the conditional event. Then,  the natural way of defining the indicator of $A|H$ (denoted by the same symbol) is (see, e.g. \cite[Section 2.2]{GiSa21IJAR})
	\[
	A|H = AH + x \no{H} = \left\{\begin{array}{ll}
		1, &\mbox{if $AH$ is true,}\\
		0, &\mbox{if $\no{A}H$ is true,}\\
		x,  &\mbox{if $\no{H}$ is true,}	
	\end{array}
	\right.
	\]
	where $x=P(A|H)$. 
	Thus, the third value of $A|H$, \emph{void}, in numerical terms  is  the conditional probability $P(A|H)$. 
	The choice of $P(A|H)$, as third value of the indicator when $H$ is false, has been considered in other works  (see e.g. \cite{coletti02,gilio90,Jeff91,Lad95,Lad96,McGe89,StJe94}).
	Notice that, as $P(AH)=P(A|H)P(H)$, denoting by $\prev(\cdot)$ the symbol of  prevision, for the random quantity $A|H$ it holds that
	\[
	\pr(A|H)= \pr(AH + x \no{H})= P(AH) + x P(\no{H}) = x[P(H)+P(\no{H})] = x = P(A|H).
	\]
	Therefore, the conditional bet proposed by de Finetti can be interpreted as a bet on the conditional event $A|H$, or on the conditional {\em if} $H$ {\em then} $A$. In this respect, the notions of  conditional bet and  bet on a conditional coincide.   Notice that, by using the indicator $A|H$, the random gain $G=sH(A-x)$ can also be  written as $G=s(A|H-x)=s(AH+x\no{H})-sx$.\\
	There is a further {\em equivalent scheme} for assessing $P(A|H)$: \\
	(Scheme 3). Let us consider the random quantity $Y=AH+y\no{H}$, where by definition $y$ is the prevision of $Y$. In order to assess $y$, in the betting framework you agree for instance to pay $y$, by receiving  1, or 0, or $y$, according to whether $AH$ is true, or $\no{A}H$ is true, or $\no{H}$  is true.  	For  the checking of coherence {\em you discard all the cases where you receive back the amount $y$ that you paid, whatever $y$ be}, that is you discard the case where $\no{H}$ is true.
When the bet is not called off, the difference $Y-A|H$ is zero,  therefore its prevision must be 0, that is $\prev(Y-A|H )=y-x=0$; therefore $y=x=P(A|H)$ and $Y=A|H$. In other words, $A|H$ can also be introduced as the random quantity $AH+y\no{H}$ , under the following conditions: $(i)$ $y=\prev(AH+y\no{H})$; $(ii)$ in order to check coherence the case where you receive back the paid amount $y$, i.e. $\no{H}$ true, is discarded. In formal terms, under conditions $(i)$ and $(ii)$, it holds that: $AH+y\no{H} = (AH+y\no{H})|H = AH|H + y\no{H}|H = A|H$.
	\subsection{Betting Scheme for conditional probability assessments}
	Given  an arbitrary family of conditional events $\K$, a  probability assessment  on the events in  $\K$ is made by specifying a real function $P : \; \mathcal{K}
	\, \rightarrow \, \mathcal{R}$.
	For each finite  subfamily 
	$\mathcal{F} = \{E_1|H_1, \ldots, E_n|H_n\}$ of $\mathcal{K}$, the restriction of $P$ to $\F$ is denoted by  $\mathcal{P} =(p_1, \ldots, p_n)$, where  $p_i = P(E_i|H_i) \, ,\;\; i =1,\ldots,n$.  
	Within the betting scheme, the random gain associated with the single assessment $p_i$ is $G_i=s_iH_i(E_i-p_i)=s_i(E_i|H_i-p_i)$, where $s_i$ is an arbitrary real number. The random gain associated with the assessment $\P$ is $G=\sum_{i=1}^n G_i = \sum_{i=1}^n s_iH_i(E_i-p_i) = \sum_{i=1}^n s_i(E_i|H_i-p_i)$, where $s_1, \ldots, s_n$ are arbitrary  real numbers. 
	The random gain  $G$ can be looked at as the difference $\sum_{i=1}^n s_i E_i|H_i - \sum_{i=1}^n s_i p_i$, that is the difference between what you receive, $\sum_{i=1}^n s_i E_i|H_i$, and what you pay, $\sum_{i=1}^n s_i p_i$.  We set $\mathcal{H}_n=H_1\vee \cdots \vee H_n$ and we denote by $G_{\mathcal{H}_n}$ the restriction of $G$ to $\mathcal{H}_n$, that is the random gain when at least a bet is not called off. 
We observe that, by expanding the second member of the equality 
	\[\Omega=(E_1H_1 \vee \no{E}_1H_1 \vee \no{H}_1) \wedge \cdots \wedge (E_nH_n \vee \no{E}_nH_n \vee \no{H}_n),
	\]
	we obtain the disjunction of  $3^n$ conjunctions $A_1 \cdots A_n$, where $A_i \in \{E_iH_i, \no{E}_iH_i, \no{H}_i\},\, i=1,\ldots,n$. By discarding the conjunctions which coincide with the impossible event $\emptyset$, the sure event $\Omega$ can be represented as a suitable disjunction $\Omega = C_0 \vee C_1  \cdots \vee C_m$, where  $C_1, \ldots, C_m$ are the constituents which logically imply $\H_n$, with $m \leq 3^n-1$.  If $\no{H}_1\cdots \no{H}_n\neq \emptyset$, we denote it by $C_0$. When $C_0$ is true, $G=g_0=0$ (all the bets are called off). Notice that, for every $h=1,\ldots,m$, the value $g_h$ of $G_{\H_n}$ associated with the constituent $C_h$ is given by $g_h=\sum_{i:C_h\subseteq E_iH_i}s_i(1-p_i)-\sum_{i:C_h\subseteq \no{E}_iH_i}p_is_i$.  We give below the definition of coherence.  	
\begin{definition}\label{DEF:COER-CONDBET}  The function $P$ defined on an arbitrary family of conditional events $\mathcal{K}$ is said to be {\em coherent}
if and only if, for every integer $n$, for every finite subfamily $\mathcal{F}=\{E_1|H_1,\ldots,E_n|H_n\}$ of $\mathcal{K}$ and for every  real numbers $s_1, \ldots, s_n$, one has:
$\min  G_{\mathcal{H}_n} \leq 0 \leq \max G_{\mathcal{H}_n}$.
\end{definition}
	The condition $\min G_{\mathcal{H}_n} \leq 0 \leq \max G_{\mathcal{H}_n}$ can be equivalently written as: $\min G_{\mathcal{H}_n} \leq 0$, or  $\max G_{\mathcal{H}_n} \geq 0$. As shown by Definition \ref{DEF:COER-CONDBET}, a probability assessment is coherent if and only if, in any finite combination of $n$ bets, after discarding the case where all the  bets are called off, it does not happen that the values $g_1,\ldots,g_m$ of the random gain  are all positive, or all negative ({\em no Dutch Book}). 
\begin{remark}\label{REM:GAINCOND}
Given a probability assessment $\P=(p_1,\ldots,p_n)$  on a finite family $\F=\{E_1|H_1,\ldots,E_n|H_n\}$,  differently from the unconditional case (see Remark \ref{REM:GAIN}), the condition  $\min  G_{\mathcal{H}_n} \leq 0 \leq \max G_{\mathcal{H}_n}$ is necessary but not sufficient for the coherence of $\P$, as shown in the next example. 
	\end{remark}
\begin{example}
Let $E_1,E_2,H_1,H_2$ be four events, with $E_1H_1=\emptyset$, $H_1\neq \emptyset$,  $H_2\neq \emptyset$, $\no{H}_1E_2H_2\neq \emptyset$. Given the probability assessment $\P=(\frac{1}{2},1)$ on $\F=\{E_1|H_1,E_2|H_2\}$, the random gain is  $G=s_1H_1(E_1-\frac12)+s_2H_2(E_2-1)$. We observe that when the constituent $\no{H}_1E_2H_2$ is true, it holds that $G_{\H_2}=s_1\cdot 0(E_1-\frac12)+s_2\cdot 1(1-1)=0$ for every $s_1,s_2$. Thus, the condition $\min G_{\H_2}\leq 0 \leq \max G_{\H_2}$ is satisfied. However, the assessment $\P$ is incoherent because the condition 
$\min G_{H_1} \leq 0\leq \max G_{H_1}$ associated with the subassessment $P(E_1|H_1)=\frac12$  is not satisfied. Indeed, as $E_1H_1=0$, it holds that $G=-\frac12 s_1$, or $G=0$, according to whether $H_1$ is true, or false, respectively. Then 
$G_{H_1}=-\frac{1}{2}s_1<0$ for every $s_1>0$.
\end{example}
\subsection{On axiomatic conditional probability and coherence}
We observe that, given a function $P$ defined on $\B\times \H$, where $\B$ and $\H$ are  arbitrary families of  events, with $H \neq \emptyset$ for every $H \in \H$,  $P$ satisfies  the properties of a    conditional probability if the following conditions hold.   
\begin{enumerate}
 \item[$(i)$]  $P(E|H)\geq 0$ and $P(H|H)=1$  for every $E\in\B$ and  $H\in \H$;
 \item[$(ii)$] $P[(E_1\vee E_2)|H]=P(E_1|H)+P(E_2|H)$, if  $E_1E_2H=\emptyset$, $E_1,E_2\in \B$, $H\in \H$;
 \item[$(iii)$] $P(EH|K) = P(E|HK)P(H|K)$,  for every $E,EH,H\in \mathcal{B}$ and  $K,HK \in \H$.
 \end{enumerate}
We recall that, if $P$ is  coherent, then the properties $(i), (ii)$, and $(iii)$ are satisfied; however, these properties are in general not sufficient for coherence of $P$. 
We list below some  conditions  which, together with 
$(i), (ii)$, and $(iii)$,
imply coherence of $P$ (see, e.g.,\cite{coletti02,Dubins75,gili89,Gili95a,GiSp95,Holz85,Rega85,Rigo88,Sanf12}): \\
$(a)$ $\mathcal{B}$ is a boolean algebra  and $\H=\B^0=\mathcal{B}\setminus\{\emptyset\}$;
\\
$(b)$ $\mathcal{B}$ is a boolean algebra  and $\H\cup \{\emptyset\}$ is a subalgebra of $\mathcal{B}$;\\
$(c)$ $\mathcal{B}$ is a boolean algebra  and $\H$ is an additive (or $P$-quasi additive) subset of $\mathcal{B}\setminus\{\emptyset\}$;\\
$(d)$ $\mathcal{B}$ is a boolean algebra, $\H$ is a nonempty subset of  $\mathcal{B}\setminus\{\emptyset\}$,  and the following condition is satisfied (\cite{Csas55}) 
\[
\prod_{i=1}^nP(E_i|H_i)=\prod_{i=1}^nP(E_i|H_{i+1}),
\]
where $E_i\in \B$, $H_i\in \H$, $E_i\subseteq H_iH_{i+1}$, $i=1,\ldots, n$, and $H_{n+1}=H_1$.
\\
When $P$ is coherent, it is called a {\em conditional probability}. 
In particular, 
when $\B$ is an algebra and $\H=\B^0$, $P$  is said a full conditional probability (\cite{Dubins75}). 
Moreover, if a function $P$ defined on an arbitrary family of conditional events $\K$ is coherent, then it  can be extended to a full  conditional probability $P'$ defined on $\B\times \B^0\supseteq \K$.
\\  \ \\

{\em Geometrical representation of coherence.}  
Let $\mathcal{P} =(p_1, \ldots, p_n)$, where  $p_i = P(E_i|H_i) \, ,\;\; i =1,\ldots,n$, a probability assessment on 
	$\mathcal{F} = \{E_1|H_1, \ldots, E_n|H_n\}$. 
We represent each constituent $C_h$ by a point $Q_h=(q_{h1},\ldots,q_{hn})$, where (for each index $i$) $q_{hi}=1$, or $q_{hi}=0$, or $q_{hi}=p_i$, according to whether $C_h \subseteq E_i$, or $C_h \subseteq \no{E}_i$, or $C_h \subseteq \no{H}_i$, respectively. Notice that $q_{hi}$ is the value assumed by the indicator $E_i|H_i$ when $C_h$ is true; then, the possible values of the random vector $(E_1|H_1, \ldots, E_n|H_n)$ are the points $Q_1, \ldots, Q_m,Q_0$, with $Q_0=\P$. The points $Q_h$'s were introduced in (\cite{gilio90}, see also \cite{Gili89C0}) with the name of $(\F,\P)$-atoms, or generalized constituents. 
For each constituent $C_h$, the associated  value $g_h$ of the random gain $G=\sum_{i=1}^n s_i (E_i|H_i - p_i)$ is given by  $g_h=\sum_{i=1}^n s_i(q_{hi}-p_i)=f(Q_h - \P)=f(Q_h)-f(\P)$, where the function $f$ is defined as $f(x_1,\ldots,x_n)=s_1x_1 + \cdots + s_nx_n$. Notice that, for each constant $k$, the equation $s_1x_1 + \cdots + s_nx_n=k$ describes an hyperplane in $\mathbb{R}^n$. In particular, $\P$ belongs to the hyperplane $\pi_{\P}$ with equation $s_1x_1 + \cdots + s_nx_n=f(\P)$; moreover, each $Q_h$ belongs to the hyperplane $\pi_h$ with equation $s_1x_1 + \cdots + s_nx_n=f(Q_h)$. We observe that the hyperplanes $\pi_{\P},\pi_1,\ldots,\pi_m$ are parallel.
We denote by $\I$ the convex hull $Q_1,\ldots, Q_m$; then, we distinguish two cases: $(i)$ $\P\notin \I $; $(ii)$ $\P\in \I $. \\
In the case $(i)$, there exists an hyperplane $\pi^*$, parallel to $\pi_{\P},\pi_1,\ldots,\pi_m$, which separates $\P$ from $\I$ and hence  
\[
f(Q_h)>k^*>f(P) \,,\; h=1,\ldots,m, \;\; \mbox{or} \;\; f(Q_h)<k^*<f(P) \,,\; h=1,\ldots,m.
\]
Then  
\[
g_h=f(Q_h)-f(P) > 0 \,,\; h=1,\ldots,m, \;\; \mbox{or} \;\; g_h=f(Q_h)-f(P) < 0 \,,\; h=1,\ldots,m,
\]
that is $\min G_{\H_n} > 0$, or $\max G_{\H_n} < 0$. Therefore, the condition $\P \notin \I$ implies that the condition of coherence $\min G_{\H_n} \leq 0 \leq \max G_{\H_n}$ is not satisfied. \\
In the case $(ii)$ there exists a vector $\Lambda=(\lambda_1,\ldots,\lambda_m)$ such that 
\[
\P = \sum_{h=1}^m \lambda_h Q_h \,,\;\; \sum_{h=1}^m \lambda_h =1 \,,\;\; \lambda_h \geq 0 \,,\; h=1,\ldots,m. 
\]
Then 
\[
\sum_{h=1}^m \lambda_h g_h = \sum_{h=1}^m \lambda_h [f(Q_h)-f(\P)] = \sum_{h=1}^m \lambda_h f(Q_h) - f(\P) = f(\P) - f(\P) = 0\,,
\]
and hence: $\min \{g_1,\ldots,g_m\} \leq 0 \leq \max \{g_1,\ldots,g_m\}$; that is, $\min \G_{\H_n} \leq 0 \leq \max \G_{\H_n}$. Therefore, $\P \in \I$ implies that $\min G_{\H_n} \leq 0 \leq \max G_{\H_n}$. In conclusion, based on the analysis of cases $(i)$ and $(ii)$, it holds that: 
\begin{equation}\label{EQ:CONV-GUAD-COND}
\P \in \I \; \Longleftrightarrow \;  \min G_{\H_n} \leq 0 \leq \max G_{\H_n}.
\end{equation}
Formula (\ref{EQ:CONV-GUAD}) also follows
by exploiting a suitable alternative theorem  (\cite[Theorem 2.9]{Gale60}, see also \cite[Theorem 1]{GiSa21IJAR}).
\begin{remark}\label{REM:EQUIV-CONV-GUAD}
    We observe that when $H_i=\Omega$, $i=1,\ldots,n$, it holds that  $\H_n=\Omega$. Moreover, all points $Q_h$'s have binary components, that is   $Q_h\in\{0,1\}^n$, $h=1,\ldots,n$. Then, in the particular case where $H_i=\Omega$, $i=1,\ldots,n$, formula 
     (\ref{EQ:CONV-GUAD-COND}) becomes formula (\ref{EQ:CONV-GUAD}).
\end{remark}
Based on Definition \ref{DEF:COER-CONDBET} and formula (\ref{EQ:CONV-GUAD-COND}) we obtain the following 
 geometrical characterization of coherence.
\begin{theorem}\label{TH:CONVHULLS}
The function $P$ defined on an arbitrary family of conditional events $\mathcal{K}$ is   coherent
if and only if, for every integer $n$, for every finite subfamily $\mathcal{F}=\{E_1|H_1,\ldots,E_n|H_n\}$ of $\mathcal{K}$, denoting by $\P$ the restriction of $P$ to $\F$, it holds that $\P\in \I$.
\end{theorem}
The next example illustrates the previous characterization. 
\begin{example}\label{EX:CONV-INC}
Given four events $A,H,B,K$, with $H \neq \emptyset, K \neq \emptyset$, assume that $AH\no{K} = AH\no{B}K = \no{H}\no{B}K = \no{A}HBK = \emptyset$. Then, the constituents generated by the family $\F=\{A|H,B|K\}$ are 
\[
C_0=\no{H}\no{K} \,,\; C_1=AHBK \,,\; C_2=\no{H}BK \,,\; C_3=\no{A}H\no{K} \,,\; C_4=\no{A}H\no{B}K \,,
\]
and, by setting $P(A|H)=x, P(B|K)=y$, the points $Q_1, Q_2, Q_3, Q_4$ are 
\[
Q_1=(1,1) \,,\; Q_2=(x,1) \,,\; Q_3=(0,y) \,,\; Q_4=(0,0) \,,
\]
with $Q_0=\P=(x,y)$. 
 First of all we observe that, concerning  the sub-assessment $x$ (resp., $y$) on $A|H$ (resp., $B|K$), the condition $\P \in\I $ reduces to $x\in[0,1]$ (resp., $y\in[0,1]$); thus $(x,y)\in[0,1]^2$.
Given any $(x,y)\in[0,1]^2$,   the condition $\P \in \I$ is satisfied if and only if $ x \leq y$. Indeed, if    $x=y+a$, with $a>0$, the line with equation $x=y+\frac{a}{2}$ separates the $\P$ from the convex hull $\I$ and hence $\P$ is not coherent. On the contrary, if $x\leq y$, then $\P\in \I$. Then,  $(x,y)$ is coherent if and only if $0\leq x\leq y\leq  1$. 
\end{example}
Notice that in Example \ref{EX:CONV-INC} it holds that  $A|H\subseteq B|K$, that is $A|H$ true implies $B|K$ true and $B|K$ false implies $A|H$  false (Goodman-Nguyen inclusion relation between two conditional events).
\subsection{Penalty Criterion for conditional probability assessments} 
	The notion of coherence can be also defined by exploiting a penalty criterion, which is equivalent to the betting scheme. With the assessment $\mathcal{P} =(p_1, \ldots, p_n)$ on $\mathcal{F}=\{E_1|H_1,\ldots,E_n|H_n\}$ we associate a random loss $\L=\sum_{i=1}^n H_i(E_i-p_i)^2=\sum_{i=1}^n (E_i|H_i-p_i)^2$, which represents the square of the distance between the random point $(E_1|H_1,\ldots,E_n|H_n)\in\{0,1,p_i\}^n$ and the prevision point $\P$. More precisely, denoting by $L_h$ the value of $\L$ associated with the constituent $C_h$, it holds that $L_h=\no{\P Q_h}^2$, with $L_0=\no{\P Q_0}^2=0$.
	Given another prevision point $\mathcal{P}^* =(p_1^*, \ldots, p_n^*)$, we denote by $\L^*$ the loss associated with $\P^*$. Then, within the  {\em penalty criterion}, we
	have
\begin{definition}\label{DEF:COER-CONDPENALTYn}   
A probability assessment $\P=(p_1,\ldots,p_n)$  on a  family $\F=\{E_1|H_1,\ldots,E_n|H_n\}$ is said to be {\em coherent} if and only if 
		there does not exist $\mathcal{P}^* =(p_1^*, \ldots, p_n^*)$ such that: 
		$L_h^* \leq L_h$, for every $h$, with $L_h^* < L_h$ for at least an index $h$.
\end{definition} 

{\em Geometrical interpretation of coherence.} 
Given a probability  assessment $\P=(p_1,\ldots,p_n)$  on a  family $\F=\{E_1|H_1,\ldots,E_n|H_n\}$, for each subfamily $\S$ of $\F$, we denote by $\P_{\S}$ the sub-assessment of $\P$ on $\S$ and by $\I_{\S}$ the convex-hull associated with the pair $(\S,\P_{\S})$. Then, based on Definition \ref{DEF:COER-CONDPENALTY}, it can be verified that \cite[Theorem 4.4]{gilio90}
\begin{equation}\label{EQ:COND-CONV-PEN}
\P \text{ is coherent  } \Longleftrightarrow\;  \P_{\S} \in\I_{\S},\; \forall \S\subseteq \F\,.
\end{equation}
The general definition of coherence with the penalty criterion is given below
	\begin{definition}\label{DEF:COER-CONDPENALTY}  A function $P$ defined on an arbitrary family of conditional events $\mathcal{K}$ is said to be {\em coherent} if and only if, for every integer $n$, for every subfamily $\mathcal{F}=\{E_1|H_1,\ldots,E_n|H_n\} \subseteq \mathcal{K}$,
		denoting by $\P=(p_1,\ldots,p_n)$  the restriction of $P$ to $\F$, 
		there does not exist $\mathcal{P}^* =(p_1^*, \ldots, p_n^*)$ such that:
		$\L^*\leq \L$ and $\L^*\neq \L$, that is
		$L_h^* \leq L_h$, for every $h$, with $L_h^* < L_h$ for at least an index $h$.
	\end{definition} 
Based on formulas   (\ref{EQ:CONV-GUAD-COND}) and (\ref{EQ:COND-CONV-PEN}), the notions of  coherence with the betting scheme and the penalty criterion are equivalent, because in both cases coherence of the function $P$ on $\K$ amounts to (\cite{Gili89C0},\cite{gilio90})
\[
\P \in \I,\;\; \text{ for every finite subfamily } \F\subseteq \K.
\]
We observe that the loss $\L$ can be written as  $\L=\sum_{i=1}^n{\L_i}$, where
$\L_i=H_i(E_i-p_i)^2$ is the loss associated with the assessment $P(E_i|H_i)=p_i$. 
A generalization of Definition \ref{DEF:COER-CONDPENALTY} is obtained if the loss $H_i(E_i-p_i)^2$ is replaced by  $H_is(E_i,p_i)$, where $s(E_i,p_i)$ is a given (strictly) proper scoring rule.
A scoring rule for the probability of a given event $E_i$ is a function of the indicator $E_i$ and of the assessed probability $P(E_i)=p_i$. Assume that you were asked to assert $P(E_i)$, knowing that your assertion were to be scored according to the rule $s(E_i, P (E_i))$; moreover, assume that your degree of belief were $P(E) = p$, while you announced instead some other number $P (E) = x$, in the expectation that you would achieve a ``better'' score. The rule is said to be proper if you cannot expect a better score by specifying a value $x$ different from $p$. A score may represent a reward or a penalty; we think of scores as penalties, so that to improve the score means to reduce it. 
The scoring rule $s(E_i,p_i)=(E_i-p_i)^2$, used in Definition \ref{DEF:COER-CONDPENALTY}, is  proper and it is called Brier quadratic scoring rule. Then, we have (\cite{GiSa11a})

\begin{definition}\label{DEF:COER-CONDPENALTYSCORE} 
Let $s$ be a proper scoring rule and  $P$ be a function  defined on an arbitrary family of conditional events $\mathcal{K}$. For every integer $n$, for every subfamily $\mathcal{F}=\{E_1|H_1,\ldots,E_n|H_n\} \subseteq \mathcal{K}$, we denote  by $\P=(p_1,\ldots,p_n)$  the restriction of $P$ to $\F$ and by $\L=\sum_{i=1}^nH_is(E_i,p_i)$ the score associated with $\P$. The function $P$ 
is said to be {\em coherent} if and only if, for every  $\mathcal{F} \subseteq \mathcal{K}$,
		there does not exist $\mathcal{P}^* =(p_1^*, \ldots, p_n^*)$ such that: 
		$\L^*\leq \L$ and $\L^*\neq \L$, where $\L^*=\sum_{i=1}^nH_is(E_i,p_i^*)$.
	\end{definition} 
We denote by $\Pi$ the set  of coherent conditional probability assessments $P$  on $K$ and by 
$\Pi_s$ the set  of coherent conditional probability  associated with the proper scoring rule $s$. 
Then, it holds that  
$
\Pi_s=\Pi,
$
for every (bounded continuous strictly) proper scoring rule $s$ (\cite[Theorem 4]{GiSa11a}).
\section{Trivalent logics and compound conditionals}
\label{SEC:TRIVALENT}
In this section we analyze some notions of conjunctions defined in the setting of trivalent logics, where the conjunction is a  conditional event with set of truth values $\{\text{true}, \text{false},\text{void}\}$ (see, e.g., \cite{ciucci2012,ciucci2013}). We  check if they satisfy some basic logic and probabilistic  properties. In particular, we  examine the following conjunctions:    Kleene-Lukasiewicz-Heyting conjunction ($\wedge_K$), or de Finetti conjunction (\cite{deFi35});  Lukasiewicz conjunction  ($\wedge_L$);  Bochvar internal conjunction, or Kleene weak conjunction ($\wedge_B$);  Sobocinski conjunction, or quasi conjunction ($\wedge_S$).

\begin{itemize}
	\item   $(A|H) \wedge_{K} (B|K)=AHBK|(AHBK\vee \overline{A}H\vee \overline{B}K)$;
   \item $(A|H)\wedge_L (B|K)=AHBK|( AHBK \vee \overline{A} \,H \vee \overline{B} \, K\vee   \overline{H} \, \overline{K})$;
	\item   $(A|H)\wedge_{B} (B|K)=AHBK|HK=AB|HK$;
	\item 	$(A|H) \wedge_{S} (B|K)=[(AH\vee\overline{H})\wedge (BK\vee \overline{K})]|(H\vee K)=$ \\ $=[(A\vee\overline{H})\wedge (B\vee \overline{K})]|(H\vee K)$.
\end{itemize}
We recall that, given any conditional event $A|H$,  the negation $\no{A|H}$ is defined as $\no{A|H}=\no{A}|H$. Based on De Morgan's law we give below the list of disjunctions associated with the previous list of conjunctions. 
\begin{itemize}
	
	\item $(A|H) \vee_{K} (B|K)$ is the negation of $(\no{A}|H) \wedge_{K}(\no{B}|K)$, which coincides with $\no{A}H\no{B}K|(\no{A}H\no{B}K\vee AH\vee BK)$; then 
	\[
	(A|H) \vee_{K} (B|K) = (AH\vee BK)|(\no{A}H\no{B}K\vee AH\vee BK)\,;
	\]
	\item $(A|H)\vee_L (B|K)$ is the negation of $(\no{A}|H)\wedge_L (\no{B}|K)$, which coincides with $\no{A}H\no{B}K|(\no{A}H\no{B}K \vee AH \vee BK\vee   \overline{H} \, \overline{K})$; then  
	\[
	(A|H) \vee_{L} (B|K) = (AH \vee BK)|(\no{A}H\no{B}K \vee AH \vee BK\vee   \overline{H} \, \overline{K})  \,;
	\]
	\item  $(A|H)\vee_{B} (B|K)$ is the negation of $(\no{A}|H)\wedge_{B} (\no{B}|K)$, which coincides with $\no{A}\no{B}HK|HK=\no{A}\no{B}|HK$; then 
	\[
	(A|H)\vee_{B} (B|K) = (A \vee B)|HK \,.
	\]
	\item $(A|H) \vee_{S} (B|K)$ is the negation of	$(\no{A}|H) \wedge_{S} (\no{B}|K)$, which coincides with $(AH\vee\overline{H})\wedge (BK\vee \overline{K})|(H\vee K)$. Then 
	\[
	(A|H) \vee_{S} (B|K) = (AH \vee BK)|(H \vee K) \,.
	\]
\end{itemize}
\subsection{Some basic logical properties}
\label{SEC:EVENTI}
We consider some basic logical properties and we check if they are satisfied by the previous notions of conjunction. We start by considering the following property:\\  
P1.  Given two events $A$ and $B$ it holds that
\[
 A \subseteq B \;
 \Longleftrightarrow \;
 \; A \wedge B=A.
\] 
We replace $A$ and $B$ by the conditional events $A|H$ and $B|K$; moreover we use the Goodman-Nguyen inclusion relation among conditional events defined as (\cite[formula (3.18)]{GoNg88})
\[
A|H \subseteq B|K \; \Longleftrightarrow \; AH \subseteq BK \; \mbox{and} \; \no{B}K \subseteq \no{A}H \,.
\]
Notice that $A|H \subseteq B|K$ amounts to  $AH\no{B}K=AH\no{K}=\no{H}BK=\emptyset$.\\
We will check  the validity of P1 
for each $\wedge \in\{
\wedge_K,\wedge_L,\wedge_B,\wedge_S\}$,
by verifying the equality
 \[
A|H \subseteq B|K \; \Longleftrightarrow \; (A|H) \wedge (B|K) = A|H.
\]
By considering the conjunction $\wedge_K$, when $A|H \subseteq B|K$, it holds that $AHBK=AH$ and $\overline{A}H\vee \overline{B}K=\overline{A}H$; then,
\[
(A|H) \wedge_{K} (B|K)=AHBK|(AHBK\vee \overline{A}H\vee \overline{B}K)=AH|(AH\vee \overline{A}H) = A|H ,
\]
and hence $A|H \subseteq B|K \; \Longrightarrow \; (A|H) \wedge_{K} (B|K) = A|H$.\\
 Conversely,  if $(A|H) \wedge_{K} (B|K) = A|H$, then $A|H \subseteq B|K$ because $AH\no{B}K=AH\no{K}=\no{H}BK=\emptyset$.
Therefore $\wedge _K$ {\em satisfies} property P1, that is (see also \cite[formula (3.17)]{GoNg88}), 
 \[
A|H \subseteq B|K \; \Longleftrightarrow \; (A|H) \wedge_{K} (B|K) = A|H.
\]
Concerning the conjunction $\wedge_L$, when 
 $A|H \subseteq B|K$,  it holds that 
\[
(A|H)\wedge_L (B|K)=AHBK|( AHBK \vee \overline{A} \,H \vee \overline{B} \, K\vee   \overline{H} \, \overline{K})=AH|(H\vee \no{H}\no{K}).
\]
We observe that, when $\no{H}\no{K}$ is true, $(A|H)\wedge_L (B|K)$ is false, while $A|H$ is void; therefore  $(A|H)\wedge_L (B|K)\neq A|H$. Thus $\wedge_L$ {\em does not satisfy} the property P1.

Concerning the conjunction $\wedge_B$, when 
 $A|H \subseteq B|K$,  it holds that 
\[
(A|H)\wedge_B (B|K)=AHBK|HK=A|HK.
\]
We observe that, when $\no{A}H\no{K}$ is true, $(A|H)\wedge_B (B|K)$ is void, while $A|H$ is false; therefore  $(A|H)\wedge_B (B|K)\neq A|H$. Thus $\wedge_B$ {\em does not satisfy} the property P1.

Concerning the conjunction $\wedge_S$, when 
 $A|H \subseteq B|K$,  it holds that 
\[
(A|H)\wedge_S (B|K)=(AH\vee \no{H}BK)|(H\vee K).
\]
We observe that, when $\no{H}BK$ is true, $(A|H)\wedge_S (B|K)$ is true,  while $A|H$ is void; therefore  $(A|H)\wedge_S (B|K)\neq A|H$. Thus $\wedge_S$ {\em does not satisfy} the property~P1. \\ \ \\
P2. We now consider the following  relation between disjunction and conjunction:
\begin{equation} 
A=AB\vee A\no{B}
\end{equation}
and, when events are replaced by conditional events, we check its validity in the previous trivalent logics. We also observe that the previous equality is related to the distributivity property, because 
\[
A=A \wedge \Omega = A \wedge (B\vee \no{B}) = AB\vee A\no{B} \,.
\]
We will also examine this aspect, when replacing $A$ and $B$ by $A|H$ and $B|K$; in each trivalent logic, associated to a pair $(\wedge,\vee)$, we will check the validity of the relations 
\begin{equation}\label{EQ:P2a}
A|H=[(A|H) \wedge (B|K)]\vee [(A|H) \wedge (\no{B}|K)],
\end{equation}
\begin{equation}\label{EQ:P2b}
A|H=(A|H) \wedge (K|K),
\end{equation}
\begin{equation}\label{EQ:P2c}
(A|H) \wedge [(B|K)\vee (\no{B}|K)] = [(A|H) \wedge (B|K)]\vee [(A|H) \wedge (\no{B}|K)]  \,.
\end{equation}

Concerning  $\wedge_K$ and $\vee_{K}$, we have to check if it holds that
\[
(A|H)= [(A|H) \wedge_{K} (B|K)]\vee_K [(A|H) \wedge_{K} (\no{B}|K)] \,.
\]
We observe that 
\[
(A|H) \wedge_{K} (\no{B}|K) = AH\no{B}K|(AH\no{B}K\vee \no{A}H\vee BK) .
\]
Then
\begin{equation}\label{EQ:DECOMPK}
\begin{array}{ll}
[(A|H) \wedge_{K} (B|K)]\vee_K [(A|H) \wedge_{K} (\no{B}|K)]=\\ \ \\
=[AHBK|(AHBK\vee \no{A}H\vee \no{B}K)] \vee_K [AH\no{B}K|(AH\no{B}K\vee \no{A}H\vee BK)] =\\
=(AHBK\vee AH\no{B}K)|(AHBK\vee AH\no{B}K\vee(\no{A}H \vee \no{B}K)\wedge (\no{A}H\vee BK) )= \\
=AHK|(AHK \vee \no{A}H)\neq A|H,
\end{array}
\end{equation}
because, for instance,  when $AH\no{K}$ is true, $AHK|(AHK \vee \no{A}H)$ is void, while $A|H$ is true. Thus, formula  (\ref{EQ:P2a}) {\em is not  satisfied} by $\wedge_K$ and $\vee_K$.\\
Concerning formula  (\ref{EQ:P2b}), we observe that {\em it 
 is not satisfied}; indeed
\[
(A|H) \wedge_K (K|K) = AHK|(AHK\vee \overline{A}H) \neq A|H \,,
\]
as shown in (\ref{EQ:DECOMPK}) (notice that $(A|H) \wedge_K (H|H) = A|H$). 
Concerning formula  (\ref{EQ:P2c}), we observe that 
\[
(A|H) \wedge_K [(B|K)\vee_K (\no{B}|K)] = (A|H) \wedge_K (\Omega|K)] = AHK|(AHK\vee \overline{A}H) \,,
\]
which, by recalling (\ref{EQ:DECOMPK}), coincides with  
\[
[(A|H) \wedge_K (B|K)]\vee_K [(A|H) \wedge_K (\no{B}|K)]  \,.
\]
Then, formula  (\ref{EQ:P2c}) {\em is satisfied} by $\wedge_K$ and $\vee_K$. \\
Concerning  $\wedge_L$ and $\vee_{L}$, we have to check if it holds that
\[
(A|H)= [(A|H) \wedge_{L} (B|K)]\vee_L [(A|H) \wedge_{L} (\no{B}|K)] \,.
\]
We observe that 
\[
(A|H) \wedge_{L} (\no{B}|K) = AH\no{B}K|(AH\no{B}K\vee \no{A}H\vee BK \vee \no{H}\no{K}) .
\]
Then, as it can be verified,
\begin{equation}\label{EQ:DECOMPL}
\begin{array}{ll}
[(A|H) \wedge_{L} (B|K)]\vee_L [(A|H) \wedge_{L} (\no{B}|K)]= \\ =[AHBK|(AHBK\vee \no{A}H\vee \no{B}K \vee \no{H}\no{K})] \vee_L [AH\no{B}K|(AH\no{B}K\vee \no{A}H\vee BK\vee \no{H}\no{K})]
 = \\ = \cdots
=AHK|(AHK \vee \no{A}H\vee \no{H}\no{K}\vee AH\no{K})=AHK|(H\vee  \no{K})\neq A|H,
\end{array}
\end{equation}
because, for instance,  when $AH\no{K}$ is true, $AHK|(H\vee  \no{K})$ is false, while $A|H$ is true. Thus, formula (\ref{EQ:P2a}) {\em is not  satisfied} by $\wedge_L$ and $\vee_L$.\\
Concerning formula  (\ref{EQ:P2b}), we observe that {\em it 
 is not satisfied}; indeed
\[
(A|H) \wedge_L (K|K) = AHK|(AHK \vee \no{A}H \vee \no{H}\no{K}) \neq A|H \,,
\]
because, when  $AH\no{K}$ is true, $A|H$ is true, while $AHK|(AHK\vee \overline{A}H \vee \no{H}\no{K})$ is void (notice that $(A|H) \wedge_L (H|H) = AH \neq A|H$; more in general $(A|H)\wedge_L(B|H)=ABH$, in particular $(A|H)\wedge_L(A|H)=AH$; thus $(A|H)\wedge_L(A|H)\neq A|H$). 
Concerning formula  (\ref{EQ:P2c}), we observe that 
\[
(A|H) \wedge_L [(B|K)\vee_L (\no{B}|K)] = (A|H) \wedge_L (\Omega|K) = AHK|(AHK\vee \overline{A}H  \vee \no{H}\no{K}) \,,
\]
which, by recalling (\ref{EQ:DECOMPL}), does not coincide with
\[
[(A|H) \wedge_L (B|K)]\vee_L [(A|H) \wedge_L (\no{B}|K)]=AHK|(H\vee \no{K}),
\]
because, for instance, when $AH\no{K}$ is true, $AHK|(H\vee \no{K})$ is false, while $AHK|(AHK\vee \overline{A}H  \vee \no{H}\no{K})$ is void.
Then, formula  (\ref{EQ:P2c}) {\em is not satisfied} by $\wedge_L$ and $\vee_L$. \\
Concerning  $\wedge_B$ and $\vee_B$, we have to check if it holds that
\[
(A|H)= [(A|H) \wedge_B (B|K)]\vee_B [(A|H) \wedge_B (\no{B}|K)] \,.
\]
We observe that 
\[
(A|H) \wedge_B (B|K) = AB|HK \,,\;\;\; (A|H) \wedge_B (\no{B}|K) = A\no{B}|HK.
\]
Then
\begin{equation}\label{EQ:DECOMPB}
\begin{array}{ll}
[(A|H) \wedge_B (B|K)]\vee_B [(A|H) \wedge_B (\no{B}|K)]
= (AB \vee A\no{B})|HK = A|HK \neq A|H \,,
\end{array}
\end{equation}
because, for instance,  when $AH\no{K}$ is true, $A|HK$ is void, while $A|H$ is true. 
Thus, formula (\ref{EQ:P2a}) {\em is not  satisfied} by $\wedge_B$ and $\vee_B$. \\
Concerning formula  (\ref{EQ:P2b}), we observe that {\em it 
 is not satisfied}; indeed
\[
(A|H) \wedge_B (K|K) = A|HK \neq A|H \,,
\]
as shown in (\ref{EQ:DECOMPB})
 (notice that $(A|H) \wedge_B (H|H) = A|H$; more in general $(A|H)\wedge_B(B|H)=AB|H$). 
Concerning formula  (\ref{EQ:P2c}), we observe that 
\[
(A|H) \wedge_B [(B|K)\vee_B (\no{B}|K)] = (A|H) \wedge_B (\Omega|K) = A|HK \,,
\]
which, by recalling (\ref{EQ:DECOMPB}),  coincides with
\[
[(A|H) \wedge_B (B|K)]\vee_B [(A|H) \wedge_B (\no{B}|K)].
\]
Then, formula  (\ref{EQ:P2c}) {\em is  satisfied} by $\wedge_B$ and $\vee_B$. \\

Concerning  $\wedge_S$ and $\vee_S$, we have to check if it holds that
\[
(A|H)= [(A|H) \wedge_S (B|K)]\vee_S [(A|H) \wedge_S (\no{B}|K)] \,.
\]
We observe that 
\[
(A|H) \wedge_S (B|K) = (AHBK \vee AH\no{K} \vee \no{H}BK)|(H \vee K) \,,
\] 
and 
\[
(A|H) \wedge_S (\no{B}|K) = (AH\no{B}K \vee AH\no{K} \vee \no{H}\no{B}K)|(H \vee K).
\]
Then
\begin{equation}
\label{EQ:DECOMPS}
\begin{array}{ll}
[(A|H) \wedge_S (B|K)]\vee_S [(A|H) \wedge_S (\no{B}|K)]
=(AH \vee \no{H}K)|(H \vee K) = \\ \ \\ = (A \vee \no{H})|(H \vee K) \neq A|H \,,
\end{array}
\end{equation}
because, for instance,  when $\no{H}K$ is true, $(A \vee \no{H})|(H \vee K)$ is true, while $A|H$ is void. Thus, formula (\ref{EQ:P2a}) {\em is not  satisfied} by $\wedge_S$ and $\vee_S$. \\
Concerning formula  (\ref{EQ:P2b}), we observe that {\em it 
 is not satisfied}; indeed
\[
(A|H) \wedge_S (K|K) = (A\vee \no{H})|(H\vee K) \neq A|H \,,
\]
as shown in (\ref{EQ:DECOMPS}).
 (notice that $(A|H) \wedge_S (H|H) = A|H$; more in general $(A|H)\wedge_S(B|H)=AB|H$). 
Concerning formula  (\ref{EQ:P2c}), we observe that 
\[
(A|H) \wedge_S [(B|K)\vee_S (\no{B}|K)] = (A|H) \wedge_S (K|K) = (A\vee\no{H})|(H\vee K) \,,
\]
which, by recalling (\ref{EQ:DECOMPS}),  coincides with
\[
[(A|H) \wedge_S (B|K)]\vee_S [(A|H) \wedge_S (\no{B}|K)].
\]
Then, formula  (\ref{EQ:P2c}) {\em is  satisfied} by $\wedge_S$ and $\vee_S$. \\
P3. We now consider the following  relation between disjunction and conjunction:
\[
A \vee B = A \vee \no{A}B=B \vee A\no{B}
\]
We observe that, in logical terms, the relation $A \vee B = A + B - AB$ is equivalent to $A \vee B = A \vee \no{A}B = B \vee A\no{B}$; we will check this kind of property. \\
Concerning the conjunction $\wedge_{K}$ and disjunction $\vee_{K}$, we have to check if it holds that
\[
(A|H) \vee_{K} (B|K) = (A|H) \vee_{K} [(\no{A}|H) \wedge_{K} (B|K)] \,.
\]
We observe that 
\[
(\no{A}|H) \wedge_{K} (B|K) = \no{A}HBK|(\no{A}HBK\vee AH\vee \overline{B}K) \,;
\]
then 
\[
(A|H) \vee_{K} [(\no{A}|H) \wedge_{K} (B|K)] = (A|H) \vee_{K} \no{A}HBK|(\no{A}HBK\vee AH\vee \overline{B}K)= 
\]
\[
= (AH \vee \no{A}HBK)|(AH \vee \no{A}HK) \,.
\]

We observe that, when $\no{H}BK$ is true, $(AH \vee \no{A}HBK)|(AH \vee \no{A}HBK \vee \no{A}H\no{B}K)$ is void, while $(A|H) \vee_{K} (B|K)$ is true; therefore 
\[
(A|H) \vee_{K} (B|K) \; \neq \; (A|H) \vee_{K} [(\no{A}|H) \wedge_{K} (B|K)] \,.
\]
By a similar reasoning, it can be verified that 
\[
(A|H) \vee_{K} (B|K) \; \neq \; (B|K) \vee_{K} [(A|H) \wedge_{K} (\no{B}|K)] \,.
\]
Indeed, 
\[
(A|H) \wedge_{K} (\no{B}|K) = AH\no{B}K|(AH\no{B}K\vee \no{A}H\vee BK) \,,
\]
and 
\[
(B|K) \vee_{K} [(A|H) \wedge_{K} (\no{B}|K)] = (BK \vee AH\no{B}K)|(BK \vee AH\no{B}K \vee \no{A}H\no{B}K) \,.
\]
Then, when $AH\no{K}$ is true, $(BK \vee AH\no{B}K)|(BK \vee AH\no{B}K \vee \no{A}H\no{B}K)$ is void, while $(A|H) \vee_{K} (B|K)$ is true. 
Thus, $\vee_{K}$ and $\wedge_{K}$ {\em do not satisfy} the property P3. \\
Concerning the conjunction $\wedge_{L}$ and disjunction $\vee_{L}$, we have to check if it holds that
\[
(A|H) \vee_{L} (B|K) = (A|H) \vee_{L} [(\no{A}|H) \wedge_{L} (B|K)] \,.
\]
We observe that 
\[(\no{A}|H)\wedge_L (B|K)=\no{A}HBK|(\no{A}HBK \vee AH \vee \overline{B} \, K\vee   \overline{H} \, \overline{K}) \,;
\]
moreover, it can be verified that
\[
(A|H) \vee_{L} [\no{A}HBK|(\no{A}HBK \vee AH \vee \overline{B} \, K\vee   \overline{H} \, \overline{K})] = (AH\vee \no{A}HBK)|( \no{A}HK \vee AH \vee \no{H}BK) \,.
\]
Then, when $\no{H}BK$ is true, $(AH\vee \no{A}HBK)|( \no{A}HK \vee AH \vee \no{H}BK)$ is false, while $(A|H) \vee_{L} (B|K)$ is true. Therefore 
\[
(A|H) \vee_{L} (B|K) \; \neq \; (A|H) \vee_{L} [(\no{A}|H) \wedge_{L} (B|K)] \,.
\]
By a similar reasoning, it can be verified that 
\[
(A|H) \vee_{L} (B|K) \; \neq \; (B|K) \vee_{L} [(A|H) \wedge_{L} (\no{B}|K)] \,.
\]
Indeed, 
\[
(A|H) \wedge_{L} (\no{B}|K) = AH\no{B}K|(AH\no{B}K\vee \no{A}H\vee BK \vee \no{H}\no{K}) \,;
\]
moreover, it can be verified that 
\[
(B|K) \vee_{L} [(A|H) \wedge_{L} (\no{B}|K)] = (BK \vee AH\no{B}K)|(H\no{B}K \vee BK \vee AH\no{K}) \,.
\]
Then, when $AH\no{K}$ is true, $(BK \vee AH\no{B}K)|(H\no{B}K \vee BK \vee AH\no{K})$ is false, while $(A|H) \vee_{L} (B|K)$ is true. 
Thus, $\vee_{L}$ and $\wedge_{L}$ {\em do not satisfy} the property P3. \\
Concerning the conjunction $\wedge_{B}$ and disjunction $\vee_{B}$, by recalling that $(A|H) \vee_{B} (B|K) = (A \vee B)|HK$, it holds that
\[
(A|H) \vee_{B} (B|K) = (A|H) \vee_{B} [(\no{A}|H) \wedge_{B} (B|K)] \,.
\]
Indeed, $(\no{A}|H)\wedge_B (B|K)=\no{A}B|HK$; moreover 
\[
(A|H) \vee_{B} [(\no{A}|H)\wedge_B (B|K)] = (A|H)\vee_B (\no{A}B|HK)  = (A \vee \no{A}B)|HK  = (A \vee B)|HK\,.
\]
By a similar reasoning, it holds that
\[
(A|H) \vee_{B} (B|K) = (B|K) \vee_{B} [(A|H) \wedge_{B} (\no{B}|K)] \,.
\]
Indeed, $(A|H) \wedge_{B} (\no{B}|K)=A\no{B}|HK$; moreover 
\[
(B|K) \vee_{B} [(A|H)\wedge_B (\no{B}|K)] = (B|K)\vee_B (A\no{B}|HK)  = (B \vee A\no{B})|HK  = (A \vee B)|HK\,.
\]
Thus, $\vee_{B}$ and $\wedge_{B}$ {\em satisfy} the property P3. \\
Concerning the conjunction $\wedge_{S}$ and disjunction $\vee_{S}$, we have to check if it holds that 
\[
	(A|H) \vee_{S} (B|K) = (A|H) \vee_S [(\no{A}|H) \wedge_{S} (B|K)] \,.
\]
We observe that 
\[
(\no{A}|H) \wedge_{S} (B|K) = [(\no{A} \vee \no{H})(B \vee \no{K})]|(H \vee K) = (\no{A}B \vee \no{A}\no{K} \vee \no{H}B)|(H \vee K) \,;
\]
moreover,  
\[
(A|H) \, \vee_S \, [(\no{A}B \vee \no{A}\no{K} \vee \no{H}B)|(H \vee K)] = (AH \vee \no{A}BH \vee \no{A}BK \vee \no{A}H\no{K} \vee \no{H}BK)|(H \vee K) \,.
\]
When $\no{A}H\no{K}$ is true, $(AH \vee \no{A}BH \vee \no{A}BK  \vee \no{A}H\no{K} \vee \no{H}BK)|(H \vee K)$ is true, while $(A|H) \vee_{S} (B|K)$ is false. Therefore 
\[
(A|H)\vee_S  (B|K)\neq (A|H) \vee_S [(\no{A}|H) \wedge_{S} (B|K)].
\]
By a similar reasoning, it can be verified that 
\[
(A|H) \vee_{S} (B|K) \; \neq \; (B|K) \vee_{S} [(A|H) \wedge_{S} (\no{B}|K)] \,.
\]
Indeed, 
\[
(A|H) \wedge_{S} (\no{B}|K) = [(A \vee \no{H})(\no{B} \vee \no{K})]|(H \vee K) = (A\no{B} \vee A\no{K} \vee \no{H}\no{B})|(H \vee K) \,;
\]
moreover,  it can be verified that
\[
(B|K) \vee_{S} [(A|H) \wedge_{S} (\no{B}|K)] = (AH\vee BK \vee \no{H}K)|(H \vee K) \,.
\]
 When $\no{H}\no{B}K$ is true, $(AH\vee BK \vee \no{H}K)|(H \vee K)$ is true, while $(A|H) \vee_{S} (B|K)$ is  false. 
Thus, $\vee_{S}$ and $\wedge_{S}$ {\em do not satisfy} the property P3. 
\subsection{Some basic probabilistic properties}
\label{SEC:EVENTIPROB}
We consider some basic probabilistic properties and we check if they are satisfied by the previous notions of conjunction. We consider the following properties, which are valid for unconditional events:
\\
P4. \; $P(AB)\leq P(A)$,\,\;\; $P(A \vee B)\geq P(A)$;\\ 
P5. \; $P(A \vee B) = P(A)+P(B)-P(AB) $; \\
P6. \; $\max \, \{P(A)+P(B)-1, 0\} \leq P(AB) \leq \min \, \{P(A),P(B)\} \,$.

We replace $A$ and $B$ by the conditional events $A|H$ and $B|K$; then we check the validity of the properties P4, P5, P6, when using the conjunctions $\wedge_K, \wedge_L, \wedge_B, \wedge_S$ and disjunctions $\vee_K, \vee_L, \vee_B, \vee_S$. \\
P4. Concerning $\wedge_K$ and $\vee_K$, we have to check if it holds that 
\begin{equation}\label{EQ:INEQK}
P[(A|H) \wedge_{K} (B|K)] \leq P(A|H) \,,\;\;\;\;\; P[(A|H) \vee_{K} (B|K)] \geq P(A|H) \,.
\end{equation}
We recall that 
\[
(A|H) \wedge_{K} (B|K)=AHBK|(AHBK\vee \overline{A}H\vee \overline{B}K) \,.
\]
We set $P(A|H)=x$, $P[AHBK|( AHBK \vee \overline{A} \,H \vee \overline{B} \, K\vee   \overline{H} \, \overline{K})]=y$, with $0 \leq x \leq 1, 0\leq y \leq 1$. We can refer to the partition with constituents: 
$C_0 = \no{H}BK \vee \no{H}\no{K}$, $C_1 = AHBK$, $C_2 = AH\no{B}K$, $C_3 = AH\no{K}\no{K}$, $C_4 = \no{A}H$,$C_5 = \no{H}\no{B}K$, and associated points: $Q_0 = \P = (x,y)$, $Q_1 = (1,1)$, $Q_2 = (1,0)$, $Q_3 = (1,y)$, $Q_4 = (0,0)$, $Q_5 = (x,0)$. The convex hull $\I$ of $Q_1, \ldots, Q_5$ is the triangle with vertices $Q_1, Q_2, Q_4$ and, in order $\P=(x,y)$ be coherent, it must be $\P \in \I$, that is $y \leq x$. Indeed, this inequality also follows from the inclusion relation  
\[
AHBK|(AHBK\vee \overline{A}H\vee \overline{B}K) \; \subseteq \; A|H \,.
\]
Moreover, as $(A|H) \vee_{K} (B|K)$ is the negation of $(\no{A}|H) \wedge_{K} (\no{B}|K)$, it holds that 
\[
P[(A|H) \vee_{K} (B|K)] = 1-P[(\no{A}|H) \wedge_{K} (\no{B}|K)] \geq 1-P(\no{A}|H) = P(A|H) \,.
\]
Thus, the property P4 {\em is satisfied} by $\wedge_K$ and $\vee_K$. \\
Concerning $\wedge_L$ and $\vee_L$, we have to check if it holds that 
\begin{equation}\label{EQ:INEQL}
P[(A|H) \wedge_L (B|K)] \leq P(A|H) \,,\;\;\;\;\; P[(A|H) \vee_L (B|K)] \geq P(A|H) \,.
\end{equation}
We recall that 
\[
(A|H) \wedge_{L} (B|K)=AHBK|(AHBK\vee \overline{A}H\vee \overline{B}K \vee \no{H}\no{K}) \,.
\]
We set $P(A|H)=x$, $P[AHBK|( AHBK \vee \overline{A} \,H \vee \overline{B} \, K  \vee \no{H}\no{K})]=y$, with $0 \leq x \leq 1, 0\leq y \leq 1$. We can refer to the partition with constituents: 
$C_0 = \no{H}BK$, $C_1 = AHBK$, $C_2 = AH\no{B}K$, $C_3 = AH\no{K}$, $C_4 = \no{A}H$,$C_5 = \no{H}\no{B}K  \vee \no{H}\no{K}$. The points $Q_h$'s are the same as for the pair $(\wedge_K, \vee_K)$ and hence $\P=(x,y)$ is coherent if and only if $y \leq x$. Indeed, also in this case the inclusion relation holds, that is  
\[
AHBK|(AHBK\vee \overline{A}H\vee \overline{B}K \vee \no{H}\no{K}) \; \subseteq \; A|H \,.
\]
Moreover 
\[
P[(A|H) \vee_{L} (B|K)] = 1-P[(\no{A}|H) \wedge_{L} (\no{B}|K)] \geq 1-P(\no{A}|H) = P(A|H) \,.
\]
Thus, the property P4 {\em is satisfied} by $\wedge_L$ and $\vee_L$. \\
Concerning $\wedge_B$ and $\vee_B$, we have to check if it holds that 
\begin{equation}\label{EQ:INEQB}
P[(A|H) \wedge_B (B|K)] \leq P(A|H) \,,\;\;\;\;\; P[(A|H) \vee_B (B|K)] \geq P(A|H) \,.
\end{equation}
We recall that $(A|H) \wedge_B (B|K)=AB|HK$ and we set $P(A|H)=x$, $P(AB|HK)=y$, with $0 \leq x \leq 1, 0\leq y \leq 1$. We can refer to the partition with constituents: 
$C_0 = \no{H}$, $C_1 = AHBK$, $C_2 = AH\no{B}K$, $C_3 = AH\no{K}$, $C_4 = \no{A}HK$,$C_5 = \no{A}H\no{K}$, and associated points: $Q_0 = \P = (x,y)$, $Q_1 = (1,1)$, $Q_2 = (1,0)$, $Q_3 = (1,y)$, $Q_4 = (0,0)$, $Q_5 = (0,y)$. The convex hull $\I$ of $Q_1, \ldots, Q_5$ is the polygon with vertices $Q_1, Q_2, Q_4, Q_5$ and the condition $\P \in \I$ is satisfied for every $(x,y) \in [0,1]^2$. In other words, the assessment $(x,y)$ on $\{A|H, AB|HK\}$ is coherent, for every $(x,y) \in [0,1]^2$; indeed $AB|HK \nsubseteq A|H$. Moreover, as $(A|H) \vee_{B} (B|K) = (A \vee B)|HK$, if we set $P(A|H)=x$ and $P[(A \vee B)|HK]=z$, we can refer to the partition with constituents:  
$C_0 = \no{H}$, $C_1 = AHK$, $C_2 = AH\no{K}$, $C_3 = \no{A}H\no{B}K$, $C_4 = \no{A}HBK$, $C_5 = \no{A}H\no{K}$, and associated points: $Q_0 = \P = (x,z)$, $Q_1 = (1,1)$, $Q_2 = (1,z)$, $Q_3 = (0,0)$, $Q_4 = (0,1)$, $Q_5 = (0,z)$. The convex hull $\I$ is the polygon with vertices $Q_1, Q_2, Q_3, Q_4$ and the condition $\P \in \I$ is satisfied for every $(x,z) \in [0,1]^2$. In other words, the assessment $(x,z)$ on $\{A|H, (A \vee B)|HK\}$ is coherent, for every $(x,z) \in [0,1]^2$; indeed $A|K \nsubseteq (A \vee B)|H$. 
Thus, the property P4 {\em is not satisfied} by $\wedge_B$ and $\vee_B$. \\
Concerning $\wedge_S$ and $\vee_S$, we have to check if it holds that 
\begin{equation}\label{EQ:INEQS}
P[(A|H) \wedge_S (B|K)] \leq P(A|H) \,,\;\;\;\;\; P[(A|H) \vee_S (B|K)] \geq P(A|H) \,.
\end{equation}
We recall that $(A|H) \wedge_S (B|K)=(AH \vee \no{H})(BK \vee \no{K})|(H \vee K)$ and we set $P(A|H)=x$, $P[(AH \vee \no{H})(BK \vee \no{K})|(H \vee K)]=y$, with $0 \leq x \leq 1, 0\leq y \leq 1$. We can refer to the partition with constituents: 
$C_0 = \no{H}\no{K}$, $C_1 = AHBK \vee AH\no{K}$, $C_2 = AH\no{B}K$, $C_3 = \no{A}H$, $C_4 = \no{H}BK$,$C_5 = \no{H}\no{B}K$, and associated points: $Q_0 = \P = (x,y)$, $Q_1 = (1,1)$, $Q_2 = (1,0)$, $Q_3 = (0,0)$, $Q_4 = (x,1)$, $Q_5 = (x,0)$. The convex hull $\I$ of $Q_1, \ldots, Q_5$ is the polygon with vertices $Q_1, Q_2, Q_3, Q_4$ and the condition $\P \in \I$ is satisfied for every $(x,y) \in [0,1]^2$. In other words, the assessment $(x,y)$ on $\{A|H, (AH \vee \no{H})(BK \vee \no{K})|(H \vee K)\}$ is coherent, for every $(x,y) \in [0,1]^2$; indeed $(AH \vee \no{H})(BK \vee \no{K})|(H \vee K) \nsubseteq A|H$. Moreover, as $(A|H) \vee_S (B|K) = (AH \vee BK)|(H \vee K)$, if we set $P(A|H)=x$ and $P[(AH \vee BK)|(H \vee K)]=z$, we can refer to the partition with constituents:  
$C_0 = \no{H}\no{K}$, $C_1 = AH$, $C_2 = \no{A}HBK$, $C_3 = \no{A}H\no{B}K \vee \no{A}H\no{K}$, $C_4 = \no{H}BK$, $C_5 = \no{H}\no{B}K$, and associate3 poi4ts: $Q_0 = \P = (x,z)$, $Q_1 = (1,1)$, $Q_2 = (1,0)$, $Q_3 = (0,0)$, $Q_4 = (x,1)$, $Q_5 = (x,0)$. The convex hull $\I$ is the polygon with vertices $Q_1, Q_2, Q_3, Q_5$ and the condition $\P \in \I$ is satisfied for every $(x,z) \in [0,1]^2$. In other words, the assessment $(x,z)$ on $\{A|H, (AH \vee BK)|(H \vee K)\}$ is coherent, for every $(x,z) \in [0,1]^2$; indeed $A|H \nsubseteq (AH \vee BK)|(H \vee K)$. 
Thus, the property P4 {\em is not satisfied} by $\wedge_S$ and $\vee_S$. \\
P5. Concerning $\wedge_K$ and $\vee_K$, we have to check if it holds that 
\begin{equation}\label{EQ:TOTK}
P[(A|H) \vee_{K} (B|K)] = P(A|H) + P(B|K) -  P[(A|H) \wedge_{K} (B|K)] \,,
\end{equation}
under logical independence of $A,H,B,K$. We set $P(A|H)=x$, $P(B|K)=y$, $P[(A|H) \wedge_{K} (B|K)]=z$, with $0 \leq x \leq 1, 0\leq y \leq 1, 0\leq z \leq 1$. Moreover, we recall (\cite[Theorem 3]{SUM2018S}) that $z$ is a coherent extension of $(x,y)$ if and only if $z' \leq z \leq z''$, with $z'=0, z'' = \min \, \{x,y\}$, and hence the assessment $(x,y,z)=(1,1,0)$ is coherent. As  a Then, from (\ref{EQ:TOTK}) it would follow 
\[
P[(A|H) \vee_{K} (B|K)] = 1 + 1 - 0 = 2 \,,
\]
which is incoherent.
 Actually, in this case $P[(A|H) \vee_{K} (B|K)]=1$. Indeed, as  $P(\no{A}|H)=P(\no{B}|K)=0$  the unique coherent extension on $(\no{A}|H)\wedge_{K}(\no{B}|K)$ is $P((\no{A}|H)\wedge_{K}(\no{B}|K))=0$. As a consequence, by De Morgan Law the unique coherent extension on $(A|H)\vee_{K}(B|K)$ is $P((A|H)\vee_{K}(B|K))=1$.
Thus, the property P5 {\em is not satisfied} by $\wedge_K, \vee_K$.

Concerning $\wedge_L$ and $\vee_L$, we have to check if it holds that 
\begin{equation}\label{EQ:TOTL}
P[(A|H) \vee_L (B|K)] = P(A|H) + P(B|K) -  P[(A|H) \wedge_L (B|K)] \,,
\end{equation}
under logical independence of $A,H,B,K$. We set $P(A|H)=x$, $P(B|K)=y$, $P[(A|H) \wedge_L (B|K)]=z$, with $0 \leq x \leq 1, 0\leq y \leq 1, 0\leq z \leq 1$. Moreover, we recall (\cite[Theorem 4]{SUM2018S}) that $z$ is a coherent extension of $(x,y)$ if and only if $z' \leq z \leq z''$, with $z'=0, z'' = \min \, \{x,y\}$, and hence the assessment $(x,y,z)=(1,1,0)$ is coherent. Then, as for the pair $(\wedge_K,\vee_k)$, the property P5 {\em is not satisfied} by $\wedge_L, \vee_L$. \\
Concerning $\wedge_B$ and $\vee_B$, we have to check if it holds that 
\begin{equation}\label{EQ:TOTB}
P[(A|H) \vee_B (B|K)] = P(A|H) + P(B|K) -  P[(A|H) \wedge_B (B|K)] \,,
\end{equation}
under logical independence of $A,H,B,K$. We set $P(A|H)=x$, $P(B|K)=y$, $P[(A|H) \wedge_B (B|K)]=z$, with $0 \leq x \leq 1, 0\leq y \leq 1, 0\leq z \leq 1$. Moreover, we recall (\cite[Theorem 5]{SUM2018S}) that $z$ is a coherent extension of $(x,y)$ for every  $z \in [0,1]$, and hence the assessment $(x,y,z)=(1,1,0)$ is coherent. Then, as for the pair $(\wedge_K,\vee_k)$, the property P5 {\em is not satisfied} by $\wedge_B, \vee_B$. \\
Concerning $\wedge_S$ and $\vee_S$, we have to check if it holds that 
\begin{equation}\label{EQ:TOTS}
P[(A|H) \vee_S (B|K)] = P(A|H) + P(B|K) -  P[(A|H) \wedge_S (B|K)] \,,
\end{equation}
under logical independence of $A,H,B,K$. We set $P(A|H)=x$, $P(B|K)=y$, $P[(A|H) \wedge_S (B|K)]=z$, with $0 \leq x \leq 1, 0\leq y \leq 1, 0\leq z \leq 1$. Moreover, we recall 
(\cite{gilio12ijar}) that $z$ is a coherent extension of $(x,y)$ if and only if $z' \leq z \leq z''$, with 
\[
z' = \max \, \{x+y-1,0\} \,,\;\;\; z'' = \left\{\begin{array}{ll}
\frac{x+y-2xy}{1-xy}, &\mbox{if $(x,y) \neq (1,1)$,}\\
1, &\mbox{if $(x,y) = (1,1)$.}\\
\end{array}
\right.
\]
Then, the assessment $P[(\no{A}|H) \wedge_S (\no{B}|K)]]=w$ is a coherent extension of the assessment $(x,y)$ on $\{A|H,B|K\}$ if and only if $\nu' \leq \nu \leq \nu''$, with 
\[
\nu' = \max \, \{1-x+1-y-1,0\} = \max \, \{1-x-y,0\} \,,
\]
and
\[
\nu'' = \left\{\begin{array}{ll}
\frac{1-x+1-y-2(1-x)(1-y)}{1-(1-x)(1-y)}, &\mbox{if $(x,y) \neq (0,0)$,}\\
1, &\mbox{if $(x,y) = (0,0)$,}\\
\end{array}
\right.  =  \left\{\begin{array}{ll}
\frac{x+y-2xy}{x+y-xy}, &\mbox{if $(x,y) \neq (0,0)$,}\\
1, &\mbox{if $(x,y) = (0,0)$.}\\
\end{array}
\right.
\]
In particular, when $x=y=\frac{2}{3}$, we obtain $\nu'=0, \nu''=\frac{1}{2}$; then, the assessment $(x,y,w)=(\frac{2}{3},\frac{2}{3},\frac{1}{3})$ is coherent; moreover, as  
\[
P[(A|H) \vee_S (B|K)] = 1-P[(\no{A}|H) \wedge_S (\no{B}|K)]] \,,
\]
it follows that  
\[
P[(A|H) \vee_S (B|K)] = 1-\frac{1}{3} = \frac{2}{3} \neq \frac{2}{3}+\frac{2}{3}-\frac{1}{3} = 1 = P(A|H) + P(B|K) - P[(A|H) \wedge_S (B|K)] \,.
\]
Thus, the property P5 {\em is not satisfied} by $\wedge_S, \vee_S$. \\
P6. Concerning $\wedge_K$, we have to check if it holds that 
\begin{equation}\label{EQ:FRHOK}
\max \, \{x+y-1, 0\} \leq z \leq \min \, \{x,y\} \,,
\end{equation}
where $x=P(A|H), y=P(B|K), z=P[(A|H) \wedge_K (B|K)]$. Based on (\cite[Theorem 3]{SUM2018S}), $z$ is a coherent extension of $(x,y)$ if and only if $z' \leq z \leq z''$, with $z'=0, z'' = \min \, \{x,y\}$. As $z' \neq \max \, \{x+y-1, 0\}$, the property P6 {\em is not satisfied} by $\wedge_K$. \\
Concerning $\wedge_L$, we have to check if it holds that 
\begin{equation}\label{EQ:FRHOL}
\max \, \{x+y-1, 0\} \leq z \leq \min \, \{x,y\} \,,
\end{equation}
where $x=P(A|H), y=P(B|K), z=P[(A|H) \wedge_L (B|K)]$. Based on (\cite[Theorem 4]{SUM2018S}), $z$ is a coherent extension of $(x,y)$ if and only if $z' \leq z \leq z''$, with $z'=0, z'' = \min \, \{x,y\}$. As $z' \neq \max \, \{x+y-1, 0\}$, the property P6 {\em is not satisfied} by $\wedge_L$. \\
Concerning $\wedge_B$, we have to check if it holds that 
\begin{equation}\label{EQ:FRHOL}
\max \, \{x+y-1, 0\} \leq z \leq \min \, \{x,y\} \,,
\end{equation}
where $x=P(A|H), y=P(B|K), z=P[(A|H) \wedge_B (B|K)]$. Based on (\cite[Theorem 5]{SUM2018S}), $z$ is a coherent extension of $(x,y)$ if and only if $z' \leq z \leq z''$, with $z'=0, z'' = 1$. As $z' \neq \max \, \{x+y-1, 0\}$ and $z'' \neq \min \, \{x,y\}$, the property P6 {\em is not satisfied} by $\wedge_B$. \\
Concerning $\wedge_S$, we have to check if it holds that 
\begin{equation}\label{EQ:FRHOL}
\max \, \{x+y-1, 0\} \leq z \leq \min \, \{x,y\} \,,
\end{equation}
where $x=P(A|H), y=P(B|K), z=P[(A|H) \wedge_S (B|K)]$. We recall (\cite{gilio12ijar}) that $z$ is a coherent extension of $(x,y)$ if and only if $z' \leq z \leq z''$, with 
\[
z' = \max \, \{x+y-1,0\} \,,\;\;\; z'' = \left\{\begin{array}{ll}
\frac{x+y-2xy}{1-xy}, &\mbox{if $(x,y) \neq (1,1)$,}\\
1, &\mbox{if $(x,y) = (1,1)$.}\\
\end{array}
\right.
\]
As $z'' \neq \min \, \{x,y\}$, the property P6 {\em is not satisfied} by $\wedge_S$.
\begin{remark}
	Notice that, given the assessment $(x,y)$ on $\{A|H,B|K\}$, the lower and upper bounds on the disjunction can be obtained by exploiting De Morgan Law and the lower and upper bounds on the conjunction of $\no{A}|H$ and $\no{B}|K$, that is 
	\[
	P[(A|H)\vee (B|K)]=1-P[(\no{A}|H)\wedge (\no{B}|K)].
	\]
In  Table \ref{TAB:TABLE2} are given the  	intervals  of coherent extensions of the assessment $(x,y)$ on $\{A|H,B|K\}$ to their conjunctions and disjunctions.  As shown in Table \ref{TAB:TABLE2} the Fréchet-Hoeffding bounds  for the disjunction (i.e., $[\max\{x,y\},\min\{x+y,1\}]$) \emph{are not satisfied} by $\vee_K, \vee_L, \vee_B, \vee_S$.
\end{remark}	
	\begin{table}[!ht]
		\centering
		\begin{tabular}{|c|c|}
			\hline
			      Logical operations        & Intervals of coherent extensions  \\ \hline
			 $ \wedge_K$ &  $[0,\; \min\{x,y\}]$  $\rule{0pt}{3ex} $\\  \hline
			 $\vee_K$ &  $[\max\{x,y\}\;,1]$ $\rule{0pt}{3ex}$    \\ \hline
			 
			 $\wedge_L$ &  $[0,\; \min\{x,y\}]$ $\rule{0pt}{3ex}$\\ \hline			 
			 $\vee_L$ &      $[\max\{x,y\},\;1]$     $\rule{0pt}{3ex}$\\ \hline		
			 $\wedge_B$ & $[0,\;1]$     $\rule{0pt}{3ex}$\\ \hline		
			 $\vee_B$ & $[0,\;1]$  $\rule{0pt}{3ex}$\\ \hline		
			 $\wedge_S$ &  $\left [\max \, \{x+y-1,0\}, \; \left\{\begin{array}{ll}
			 	\frac{x+y-2xy}{1-xy}, &\mbox{if $(x,y) \neq (1,1)$}\\
			 	1, &\mbox{if $(x,y) = (1,1)$}\\
			 \end{array}
		 \right.
		 \right]$ $\rule{0pt}{5ex}$ \\\hline
			 $\vee_S$ &  $\left [ \left\{\begin{array}{ll}
	\frac{xy}{x+y-xy}, &\mbox{if $(x,y) \neq (0,0)$}\\
	0, &\mbox{if $(x,y) = (0,0)$}\\
\end{array}
\right.,\; \min \, \{x+y,1\} 
\right]$ $\rule{0pt}{5ex}$\\\hline
			 $\wedge_{gs}$ &  $ [\max \, \{x+y-1,0\}, \min\{x,y\}\;]$ $\rule{0pt}{3ex}$\\\hline
			 $\vee_{gs}$ &  $ [\max \, \{x,y\}, \min\{x+y,1\}\;]$ $\rule{0pt}{3ex}$  \\
 \hline 
		\end{tabular}
	\caption{ Intervals  of coherent extensions of the assessment $(x,y)$ on $\{A|H,B|K\}$ to their conjunctions $\wedge_K$, $\wedge_L$, $\wedge_B$, $\wedge_S$,$\wedge_{gs}$ and  their disjunctions $\vee_K$, $\vee_L$, $\vee_B$, $\vee_S$, $\vee_{gs}$. The symbols $\wedge_{gs}$ and $\vee_{gs}$ refer to conjunction and disjunction illustrated in Section \ref{SEC:CONDRAND}.}	
			\label{TAB:TABLE2}
	\end{table}	
In Table \ref{TAB:TABLE3}  we list the basic properties from P1 to P6 and, for each pair of logical operations, we insert the symbol $*$ if the property is satisfied by the pair.
\begin{table}[!ht] 
	\centering
	\begin{tabular}{|l|l|c|c|c|c|c|}
		\hline
		 &Properties        &   $\begin{array}{l}
		 	\wedge_K \\
		 	\vee_K
		 \end{array}$ &   $\begin{array}{l}\wedge_L\\ \vee_L\end{array}$ &   $\begin{array}{l}\wedge_B\\ \vee_B\end{array}$ &   $\begin{array}{l}\wedge_S\\ \vee_S\end{array}$ &$\begin{array}{l}\wedge_{gs}\\ \vee_{gs}\end{array}$\\ 
		\hline
		P1 & $\C_1\subseteq \C_2 \Longleftrightarrow \C_1\wedge\C_2=\C_1 $ & $*$ & & & & *\\
	\hline	
	 & $\C_1 = (\C_1 \wedge \C_2) \vee   (\C_1 \wedge \C_{\no{2}})$ & & & && *\\
	 	
		P2 & $\C=\C \wedge (K|K) $ & & && & *\\
		&$
			\C_1 \wedge (\C_2\vee \C_{\no{2}})=(\C_1 \wedge \C_2) \vee   (\C_1 \wedge \C_{\no{2}})
		$ & $*$ & & $*$ & $*$ &*
		\\ 
				\hline 
		P3	&	$
\C_1 \vee \C_2 =\C_1  \vee (\C_{\no{1}} \wedge \C_2) =\C_2  \vee (\C_{1} \wedge \C_{\no{2}}) \,
$ & & & $*$ & &*\\
	\hline 
P4 & $z \leq x\leq w$, $z \leq y\leq w$ & $*$ & $*$ & & 	&*	\\
		\hline
P5 & $w= x+y-z$ & & & & &	*\\
		\hline
P6 & $\max\{x+y-1,0\}\leq z\leq \min\{x,y\}$ & & & & &	*\\
 & $\max\{x,y\}\leq w\leq \min\{x+y,1\}$ & & & & &	*\\

\hline 
\end{tabular} 
\caption{Properties P1---P6 and  pairs of logical operations. 
The symbol $*$ means that the pair satisfies the property.
The symbols $\C,\C_1$, and $\C_2$ denote conditional events.
 $P(\C_1)=x$,  $P(\C_2)=y$, $P(\C_1\wedge \C_2)=z$,  $P(\C_1\vee \C_2)=w$. 
 For the pair of logical operations $(\wedge_{gs}, \vee_{gs})$,  illustrated in Section \ref{SEC:CONDRAND}, the symbols $z$ and $w$ are the \emph{previsions} of $\C_1\wedge_{gs}\C_2$ and $\C_1\vee_{gs}\C_2$, respectively.
} 	
\label{TAB:TABLE3}
\end{table}
 
\section{The conjunction and disjunction in the setting of conditional random quantities}
\label{SEC:CONDRAND}
In this section we illustrate the notions of conjunction $\wedge_{gs}$ and disjunction $\vee_{gs}$, here simply denoted by $\wedge$ and $\vee$, introduced in the setting of coherence by Gilio and Sanfilippo (see, e.g., \cite{GiSa14,GiSa19,GiSa20,GiSa21A}). This approach can be related to some works by McGee (\cite{McGe89}) and Kaufman(\cite{Kauf09}).  First of all, we observe that in probabilistic reasoning, as well as in real applications, two main aspects are: (i) the checking of consistency, or coherence, of probability assessments; (ii) the coherent extension of given initial assessments to further conditional events and/or random quantities. In both these aspects we need to suitably represent in numerical terms the truth values of conditional events. As shown in previous sections, for a conditional event $A|H$ the natural representation of the logical values true, false, void is 1, 0, $P(A|H)$, respectively; this  allows to  implement useful procedures and algorithms for checking coherence. In addition, by representing the constituents $C_h$'s by the points $Q_h$'s, we can develop a geometrical approach to coherence for both the betting scheme and the penalty criterion, as made by de Finetti in the case of unconditional events and random quantities. Based on this motivations, in our approach the notion of conjunction (as well as of disjunction) of two conditional events $A|H$ and $B|K$, $(A|H) \wedge (B|K)$, is directly defined as a suitable five-valued conditional random quantity. This numerical approach allows to distinguish three different levels of indeterminacy: $(a) \; H$ true and $K$ false; $(b) \; H$ false and $K$ true; $(c) \; H$ and $K$ both false. More precisely, the value of $(A|H) \wedge (B|K)$ is 1, or 0, or $P(A|H)$, or $P(B|K)$, or $\pr[(A|H) \wedge (B|K)]$ (the prevision of the conjunction), according to whether $AHBK$ is true, or $\no{A}H \vee \no{B}K$ is true, or $\no{H}BK$, or $AH\no{K}$ is true, or $\no{H}\no{K}$ is true, respectively. Within the betting scheme, by starting with a coherent assessment $(x,y)$ on $\{A|H, B|K\}$, if you extend $(x,y)$ (in a coherent way) by adding the assessment $\pr[(A|H) \wedge (B|K)]=z$, then you agree for instance to pay $z$, by receiving the random amount 
\[
(A|H) \wedge (B|K) = AHBK + x \no{H}BK + y AH\no{K} + z \no{H}\no{K} \; \in \; \{1,0,x,y,z\} \,.
\]
Notice that when $\no{H}\no{K}$ is true you receive back the paid amount $z$, then this case must be discarded when checking the coherence of $z$. Moreover, by recalling the reasoning in Section \ref{THIRD} on Scheme 3, the conjunction is a conditional random quantity, as written below. 
\begin{equation}\label{EQ:CONJ}
(A|H) \wedge (B|K) = (AHBK + x \no{H}BK + y AH\no{K})|(H \vee K) \,,
\end{equation} 
with 
\begin{equation}\label{EQ:MCGEERQ}
\begin{array}{ll}
\pr[(A|H) \wedge (B|K)]=\pr[(AHBK + x \no{H}BK + y AH\no{K})|(H \vee K)]=\\
=P[AHBK|(H \vee K)]+x P[\no{H}BK|(H \vee K)]+
y P[AH\no{K}|(H \vee K)].
\end{array}
\end{equation}
More in general, in order to check coherence of any given assessment $(x,y,z)$ on $\{A|H, B|K, (A|H) \wedge (B|K)\}$, we could exploit  a suitable version of Theorem \ref{TH:CONVHULLS} for conditional random quantities (see, e.g., \cite[Theorem 1]{GiSa20}). \\
Notice that, when $P(H\vee K)>0$, formula (\ref{EQ:MCGEERQ}) becomes  
\begin{equation}\label{EQ:MCGEE}
\pr[(A|H) \wedge (B|K)]=\frac{P(AHBK)+P(A|H) P(\no{H}BK)+P(B|K) P(AH\no{K})}{P(H\vee K)},
\end{equation}
which is the formula given in \cite{Kauf09,McGe89}.
In addition, if $A,B,H,K$ are stochastically independent, it holds that  
\[
\begin{array}{ll}
P(AHBK)+P(A|H) P(\no{H}BK)+P(B|K) P(AH\no{K})=\\=
P(A)P(B)P(H)P(K)+P(A)P(B)P(\no{H})P(K)+P(A)P(B)P(H)P(\no{K})=\\
=P(A)P(B)P(H\vee K),
\end{array}
\]
and hence formula (\ref{EQ:MCGEE}) reduces to
$\pr[(A|H) \wedge (B|K)]=P(A)P(B)$.\\
By a dual approach, the disjunction $(A|H) \vee (B|K)$ is defined as a conditional random quantity, with values 1, or 0, or $P(A|H)$, or $P(B|K)$, or $\pr[(A|H) \vee (B|K)]$ (the prevision of the disjunction), according to whether $AH \vee BK$ is true, or $\no{A}H\no{B}K$ is true, or $\no{H}\no{B}K$, or $\no{A}H\no{K}$ is true, or $\no{H}\no{K}$ is true, respectively. Then 
\[
(A|H) \vee (B|K) = AH \vee BK + x \no{H}\no{B}K + y \no{A}H\no{K} + w \no{H}\no{K} \; \in \; \{1,0,x,y,w\} \,,
\]
where $x=P(A|H), y=P(B|K), w=\pr[(A|H) \vee (B|K)]$. 
As for conjunction, the disjunction is a conditional random quantity, as written below. 
\begin{equation}\label{EQ:DISJ}
(A|H) \vee (B|K) = (AH \vee BK + x \no{H}\no{B}K + y \no{A}H\no{K})|(H \vee K) \,,
\end{equation} 
with $w=\pr[(A|H) \vee (B|K)]=\pr[(AH \vee BK + x \no{H}\no{B}K + y \no{A}H\no{K})|(H \vee K)]$.
It can be verified that (a numerical version of) De Morgan's laws are satisfied, that is (\cite[Theorem 5]{GiSa19})
\begin{equation}\label{EQ:DEMORGAN}
(A|H) \vee (B|K) = 1 - (\no{A}|H) \wedge (\no{B}|K) ,\,\, 
(A|H) \wedge (B|K) = 1 - (\no{A}|H) \vee (\no{B}|K) . 
\end{equation}
Now, we will show that the conjunction and disjunction, as defined in our approach, satisfy (in numerical terms) all the properties from P1 to P6. Then, we will illustrate a list of further basic properties, valid for unconditional events, which continue to hold in our approach. \\
P1.  The property $A\subseteq B \Longleftrightarrow A\wedge B=B$ in numerical terms becomes
\[
A\leq B \;\;\Longleftrightarrow \;\;A\wedge B=B.
\]
When replacing $A$  and $B$ by $A|H$ and $B|K$ it holds that (\cite[Equation (16)]{GiSa21A})
\[
A|H \leq B|K \;\; \Longleftrightarrow \;\; (A|H) \wedge (B|K) = A|H \,.
\]
Thus, the property P1 is \emph{satisfied} by $\wedge$. 
Notice that  (\cite[Equation (15)]{GiSa21A})
\[
A|H \leq B|K \;\; \Longleftrightarrow \;\; A|H \subseteq B|K \,,\; \mbox{or} \; AH=\emptyset \,,\; \mbox{or} \; BK=K \,,
\]
then $A|H \subseteq B|K \Rightarrow (A|H) \wedge (B|K) = A|H$.
%
\\
\noindent P2. We will verify that
\begin{equation}\label{EQ:P2bCOND}
A|H=(A|H) \wedge (K|K),
\end{equation}
and
\begin{equation}\label{EQ:P2cCOND}
(A|H) \wedge [(B|K)\vee (\no{B}|K)] = [(A|H) \wedge (B|K)] + [(A|H) \wedge (\no{B}|K)]  \,.
\end{equation}
Concerning formula (\ref{EQ:P2bCOND}), from (\ref{EQ:CONJ}), as $y=P(K|K)=1$,  it holds that 
\[
(A|H) \wedge (K|K) = (AHK + x \no{H}K + AH\no{K})|(H \vee K)= (AH + x \no{H}K)|(H \vee K) \,,
\]
where $x=P(A|H)$. We observe that $(AH + x \no{H}K)|(H \vee K) = A|H$ when $H \vee K$ is true,  with $(AH + x \no{H}K)|(H \vee K) = \pr[(AH + x \no{H}K)|(H \vee K)]$ when $H \vee K$ is false. Moreover,
\[
\pr[(AH + x \no{H}K)|(H \vee K)] = P[AH|(H \vee K)] + x P[\no{H}K|(H \vee K)] =
\]
\[
= x [P[H|(H\vee K)] + P[\no{H}K|(H \vee K)] = x P[(H \vee K) |(H \vee K)] = x \,;
\]
thus: $(AH + x \no{H}K)|(H \vee K) $ coincides with $A|H$ and hence formula (\ref{EQ:P2bCOND}) is satisfied. 
Concerning formula (\ref{EQ:P2cCOND}), we observe that
\[
(A|H) \wedge [(B|K)\vee (\no{B}|K)] = (A|H) \wedge (\Omega|K) = (A|H) \wedge (K|K) = A|H \,.
\]
Moreover, as it can be verified (\cite[Theorem 7 ]{GiSa20}),
\[
(A|H) \wedge (B|K) + (A|H) \wedge (\no{B}|K) = A|H ,
\]
and hence formula (\ref{EQ:P2cCOND})  is satisfied. Thus, the property P2 \emph{is satisfied}.
Notice that in (\cite{FGGS22}) a theory has been proposed on compound conditionals, where formula (\ref{EQ:P2cCOND}) can be written as 
\[
(A|H) \wedge [(B|K)\vee (\no{B}|K)] = [(A|H) \wedge (B|K)] \vee [(A|H) \wedge (\no{B}|K)]  \,.
\]
P3. In numerical terms the property $A\vee B=A\vee \no{A}B$ becomes $A\vee B=A+\no{A}B$. When replacing $A$ and $B$  by $A|H$ and $B|K$ we can verify that
\begin{equation}\label{EQ:P3-NUM}
(A|H) \vee (B|K) = (A|H) + [(\no{A}|H) \wedge (B|K)] \,.
\end{equation}
We observe that $(B|K)=(A|H) \wedge (B|K) + (\no{A}|H) \wedge (B|K)$; moreover, 
we recall that (\cite[Section 6]{GiSa14})
\begin{equation}\label{EQ:INCESCL}
	(A|H) \vee (B|K)=A|H+B|K-(A|H) \wedge  (B|K).
\end{equation}
Then,
\[
\begin{array}{l}
	(A|H) \vee (B|K) = A|H+B|K-(A|H) \wedge  (B|K)
	=(A|H)+(\no{A}|H) \wedge  (B|K).
\end{array}
\]
Thus, formula (\ref{EQ:P3-NUM}) holds, that is  the property P3 \emph{is satisfied}.\\
In the approach proposed in  (\cite{FGGS22})  the previous formula becomes 
\[
(A|H) \vee (B|K) = (A|H) \vee [(\no{A}|H) \wedge (B|K)] \,.
\]
P4.  We have to check that $ \prev[(A|H) \wedge (B|K)]\leq P(A|H)\leq \prev[ (A|H) \vee (B|K)]$.
As shown in Table \ref{TAB:TABLE1}, 
\[
 (A|H) \wedge (B|K)- A|H=-AH\no{B}K-(1-y)AH\no{K}-x\no{H}\no{B}K+(z-x)\no{H}\no{K},
\]
 where  $x=P(A|H)$, $y=P(B|K)$,  $z=\prev[(A|H) \wedge (B|K)]$. Notice that $z-x= \prev[(A|H) \wedge (B|K)- A|H]$; then,  by recalling the reasoning on  Scheme 3 in Section \ref{THIRD}, the difference $(A|H) \wedge (B|K)- A|H$ is the following conditional random quantity
 \[
 (A|H) \wedge (B|K)- A|H=(-AH\no{B}K-(1-y)AH\no{K}-x\no{H}\no{B}K)|(H\vee K),
 \]
 where $(-AH\no{B}K-(1-y)AH\no{K}-x\no{H}\no{B}K)\leq 0$.
Thus $z-x\leq 0$, that is $z\leq x$.
Likewise, $w-x\geq 0$, that is $w\geq x$, where $w=\prev[ (A|H) \vee (B|K)]$. Thus, the property P4 \emph{is satisfied} by $\wedge$ and $\vee$.
Notice that the previous results imply that
\[
(A|H) \wedge (B|K)\leq A|H\leq (A|H) \vee  (B|K).
\] 
Of course, the same inequalities hold for $B|K$.
\begin{table}[!ht]
	\centering
	\begin{tabular}{|l|l|c|c|c|c|}
		\hline
		       & $C_h$                  & $A|H$ & $B|K$ & $(A|H)\wedge(B|K)$ & $(A|H)\vee(B|K)$  \\ \hline
		$C_1 $ & $AHBK  $               &   1   &   1   &         1          &        1           \\
		$C_2 $ & $	{AH}{\no{B}K}$       &   1   &   0   &         0          &        1           \\
		$C_3 $ & $	AH\no{K}  $          &   1   &  $y$  &        $y$         &       $1$          \\
		$C_4 $ & $	{\no{A}H}{BK}$       &   0   &   1   &         0          &        1           \\
		$C_5 $ & $	{\no{A}H\no{B}K} $   &   0   &   0   &         0          &        0           \\
		$C_6 $ & $	{\no{A}H} \, \no{K}$ &   0   &  $y$  &         0          &       $y$          \\
		$C_7 $ & $	\no{H}{BK}$          &  $x$  &   1   &        $x$         &        1          \\
		$C_8 $ & $	\no{H}\,{\no{B}K}$   &  $x$  &   0   &         0          &       $x$          \\
		$C_0 $ & $	\no{H}\,\no{K} $     &  $x$  &  $y$  &        $z$         &  $w$     \\ \hline
	\end{tabular}
	\vspace{0.3cm}	
	\caption{Numerical values of $A|H$, $B|K$, $(A|H)\wedge(B|K)$, $(A|H)\vee (B|K)$  associated with the constituents $C_1,\ldots,C_8,C_0$, where $x=P(A|H)$, $y=P(B|K)$, $z=\prev[(A|H)\wedge(B|K)]$ and $w=\prev[(A|H)\vee(B|K)]$.
	}
	\label{TAB:TABLE1}
\end{table}\\
P5. From (\ref{EQ:INCESCL}) we obtain the following \emph{prevision sum rule} 
\begin{equation}\label{EQ:PINCLESCL}
\prev[(A|H) \vee (B|K)]=P(A|H)+P(B|K)-\prev[(A|H) \wedge  (B|K)].
\end{equation}
Thus, the property P5 \emph{is satisfied}  by $\wedge$ and $\vee$. Notice that the previous formula is a particular case of the more general inclusion-exclusion formula  satisfied by $\wedge$ and $\vee$ (\cite{GiSa20}). \\
P6. In \cite{GiSa14} it has been shown that, under logically independence, the set $\Pi$ of all coherent assessments $(x,y,z)$ on $\{A|H,B|K, (A|H) \wedge (B|K)\}$ is  $\{(x,y,z):(x,y)\in[0,1], \max \, \{x+y-1, 0\} \leq z \leq \min \, \{x,y\}\}$. Then, the Fr\'echet-Hoeffding bounds are satisfied when $A$ and $B$ are replaced by $A|H$ and $B|K$, that is
\[
\max \, \{x+y-1, 0\} \leq z \leq \min \, \{x,y\}.
\] 
 Thus, the property P6 \emph{is satisfied}  by $\wedge$. Moreover, by applying  (\ref{EQ:PINCLESCL}), we also obtain the dual result for the prevision $w$  of the disjunction $(A|H) \vee (B|K)$
\[
\max\{x,y\}\leq w\ \leq \min\{x+y,1\}.
\] 

\section{Some further logical and probabilistic properties}
\label{SEC:FURTHER}
We consider below some further properties which are satisfied in our approach.\\
$(a)$ We can verify that 
\begin{equation}\label{EQ:COMPTHM}
 (E|HK)\wedge (H|K)=EH|K.
\end{equation}
We  set $P(E|HK)=x$, $P(H|K)=y$,  $\prev[(E|HK)\wedge (H|K)]=z$, then  from (\ref{EQ:CONJ}) it follows that 
\[
(E|HK)\wedge (H|K)=
(EHK+x(\no{H}\vee \no{K})HK+y\no{K}EHK)|K=EH|K.
\]
Therefore property $(a)$ is satisfied;  moreover, as by  coherence $P(EH|K)=P(E|HK)P(H|K)$, it follows that 
\[
z= \prev[(E|HK)\wedge (H|K)]=P(EH|K)=P(E|HK)P(H|K)=xy.
\]
\begin{remark}
It can be verified that 
	\[
\begin{array}{ll}	
	(E|HK)\wedge_K (H|K)=EH|K,\\
	(E|HK)\wedge_L (H|K)=EHK \neq EH|K,\\
	(E|HK)\wedge_B (H|K)=E|HK \neq EH|K,\\
	(E|HK)\wedge_S (H|K)=EH|K.
\end{array}	
	\]
Therefore property $(a)$ is satisfied by $\wedge_K$ and  $\wedge_S$, while it is not satisfiedy by $\wedge_L$ and  $\wedge_B$.	 Notice that the  property $(a)$ can be equivalently written as
\begin{equation}\label{EQ:3AXIOM}
(A|B)\wedge (B|C)=A|C,  \;\; \text{ if } A\subseteq B \subseteq C.
\end{equation}
Indeed, the equality  $(E|HK)\wedge (H|K)=EH|K$,can be written as $(EHK|HK)\wedge (HK|K)=EHK|K$, with $EHK\subseteq HK \subseteq K$. Conversely, under the assumption that 
$ A\subseteq B \subseteq C$, as  $A=AB$,  $B=BC$ the equality 
$(A|B)\wedge (B|C)=A|C$ can be written as 
 $AB|C=(A|BC)\wedge (B|C)$.  \\
We recall that the equality  (\ref{EQ:3AXIOM}) is the axiom C5 in the Boolean algebra of conditionals proposed in \cite{FlGH20}.
\end{remark}	
$(b)$ By exploiting the notion of conjunction of $n$ conditional events, as defined in our approach (see, e.g., \cite{GiSa20}), we can verify that 
\begin{equation}
\begin{array}{l}
E_1\cdots E_n=E_1\wedge (E_2|E_1)\wedge \cdots \wedge (E_n|E_1\cdots E_{n-1})
\end{array} 
\end{equation} and 
\begin{equation}\label{EQ:PFACT}
	\prev[E_1\wedge (E_2|E_1)\wedge \cdots \wedge (E_n|E_1\cdots E_{n-1})]=P(E_1)P(E_2|E_1)\cdots P(E_n|E_1\cdots E_{n-1}).
\end{equation}
 The operation of  conjunction  satisfies the associativity and commutativity properties. The conjunction $\C_{1\cdots n}$ of $E_1|H_1, \ldots, E_n|H_n$ is defined as a (conditional) random quantity, with values $1$, or $0$, or $x_S$, according to whether, $\bigwedge_{i=1}^nE_iH_i$ is true, or $\bigvee_{i=1}^n\no{E}_iH_i$ is true, or $(\bigwedge_{i\notin S}E_iH_i)\wedge (\bigwedge_{i\in S}\no{H}_i)$, where $\emptyset \neq S \subseteq \{1,\ldots,n\}$ and $x_S=\prev[\bigwedge_{i\in S}(E_i|H_i)]$. In particular, when $S = \{1,\ldots,n\}$, it holds that $x_{\{1,\ldots,n\}}=\pr(\C_{1\cdots n})$, which we also denote by $x_{1\cdots n}$. In the framework of betting scheme, for instance, $x_{1\cdots n}$ is the amount to be paid in order to receive the value assumed by the random quantity $\C_{1\cdots n}$. When $\no{H}_1 \cdots \no{H}_n$ is true (i.e., $\bigvee_{i=1}^nH_i$ is false), the conditional events are all void and the value of $\C_{1\cdots n}$ coincides with the paid amount $x_{1\cdots n}$; then, for checking coherence, the case $\no{H}_1 \cdots \no{H}_n$ must be discarded. As a consequence, by also recalling the Scheme 3 in Section \ref{THIRD}, $\C_{1\cdots n}$ is a conditional random quantity with conditioning event $\bigvee_{i=1}^nH_i$.  
By a dual approach 
the disjunction $\D_{1\cdots n}$ of $E_1|H_1, \ldots, E_n|H_n$ is defined as a conditional  random quantity (with conditioning event $\bigvee_{i=1}^nH_i$), with values $1$, or $0$, or $y_S$, according to whether, $\bigvee_{i=1}^nE_iH_i$ is true, or $\bigwedge_{i=1}^n\no{E}_iH_i$ is true, or $(\bigwedge_{i\notin S}\no{E}_iH_i)\wedge (\bigwedge_{i\in S}\no{H}_i)$, where $\emptyset \neq S \subseteq \{1,\ldots,n\}$ and $y_S=\prev[\bigvee_{i\in S}(E_i|H_i)]$. In particular, when $S = \{1,\ldots,n\}$, it holds that $y_{\{1,\ldots,n\}}=\pr(\D_{1\cdots n})$, which we also denote by $y_{1\cdots n}$.  By iteratively applying (\ref{EQ:COMPTHM}) we obtain
\[
\bigwedge_{i=1}^nE_i=(\bigwedge_{i=1}^{n-1}E_i)\wedge (E_n|\bigwedge_{i=1}^{n-1}E_i)=\cdots=E_1\wedge (E_2|E_1)\wedge \cdots \wedge (E_n|E_1\cdots E_{n-1}).
\]
In addition, by recalling that
\[
P(E_1E_2\cdots E_n)=P(E_1)P(E_2|E_1)\cdots P(E_n|E_1\cdots E_{n-1}),
\]
it follows that
\[
\prev[E_1\wedge (E_2|E_1)\wedge \cdots \wedge (E_n|E_1\cdots E_{n-1})]=P(E_1)P(E_2|E_1)\cdots P(E_n|E_1\cdots E_{n-1}).
\]
The equality (\ref{EQ:PFACT})  shows that the prevision of  the conjunction of  $E_1$,  $E_2|E_1$, $\ldots$, $E_n|E_1\cdots E_{n-1}$ is equal to the product of their probabilities.\\
$(c)$
We recall that, given $n$ logically independent events $E_1,\ldots, E_n$, the best bounds on the probability of their conjunction are the
 Fr\'echet-Hoeffding bounds, that is
\begin{equation}\label{EQ:FHEV}
\max\{\sum_{i=1}^nP(E_i)-n+1,0\}\leq P(E_1\cdots E_n)\leq \min\{P(E_1),\ldots,P(E_n)\},
\end{equation}
and 
\begin{equation}\label{EQ:FHUN}
\max\{P(E_1),\ldots,P(E_n)\}\leq P(E_1\vee \cdots\vee  E_n)\leq \min\{\sum_{i=1}^nP(E_i),1\},
\end{equation}
These results still hold when replacing each event $E_i$ by the conditional event $E_i|H_i$, under logical independence of all events. In other words, in our approach the best bounds on the previsions of the conjunction and the disjunction of $n$ conditional events $E_1|H_1$, \ldots $E_n|H_n$ are
the
Fr\'echet-Hoeffding bounds (\cite{GiSa21}):
\begin{equation}\label{EQ:FHBOUNDS}
\max\{\sum_{i=1}^nP(E_i|H_i)-n+1,0\}\leq \prev[\bigwedge_{i=1}^n (E_i|H_i)]\leq \min\{P(E_1|H_1),\ldots,P(E_n|H_n)\},
\end{equation}
and
\begin{equation}\label{EQ:FHBOUNDSUN}
	\max\{P(E_1|H_1),\ldots,P(E_n|H_n)\}\leq \prev[\bigvee_{i=1}^n (E_i|H_i)]\leq  \min\{\sum_{i=1}^nP(E_i|H_i),1\}.
\end{equation}

(d)  Given $n$ events $E_1,\ldots, E_n$, we recall that by the inclusion-exclusion principle their disjunction can be represented as
\[
E_1\vee \cdots \vee E_n=\sum_{i=1}^n{E_i}-\sum_{i_1<i_2}{E_{i_1}E_{i_2}}+\sum_{i_1<i_2<i_3}E_{i_1}E_{i_2}E_{i_3}-
\cdots+(-1)^{n+1}E_1\cdots E_n.
\]
This principle still holds in our approach. Given $n$ conditional events $E_1|H_1, \ldots, E_n|H_n$, it holds that (\cite{GiSa20})
\begin{equation}\label{EQ:INCLESCLn}
\bigvee_{i=1}^n(E_i|H_i)=\sum_{i=1}^n{(E_i|H_i)}-\sum_{i_1<i_2}{(E_{i_1}|H_{i_1})\wedge (E_{i_2}|H_{i_2})}+
\cdots+(-1)^{n+1}\bigwedge_{i=1}^n(E_i|H_i).
\end{equation}
We set $P(E_i|H_i)=x_i$, $\prev[(E_{i_1}|H_{i_1})\wedge (E_{i_2}|H_{i_2})]=x_{i_1i_2}$, $i_1\neq i_2$, $\ldots$,  $\prev[\bigwedge_{i=1}^n(E_i|H_i)]=x_{1\cdots n}$, $\prev[\bigvee_{i=1}^n(E_i|H_i)]=y_{1\cdots n}$. Then, 
 when $\no{H}_1\no{H}_2\cdots \no{H}_n$ is true, that is all conditional events $E_i|H_i$ are void, the value of the left side member in (\ref{EQ:INCLESCLn}) is $y_{1\cdots n}$, while the value of the right side member is $\sum_{i=1}^nx_i-\sum_{i_1<i_2}x_{i_1i_2}+
 \cdots+(-1)^{n+1}x_{1\cdots n}$; thus the prevision of $(E_1|H_1)\vee \cdots \vee (E_n|H_n)$ is given by 
\[
y_{1\cdots n}=\sum_{i=1}^nx_i-\sum_{i_1<i_2}x_{i_1i_2}+
\cdots+(-1)^{n+1}x_{1\cdots n}.
\]
(e) Given $n$  events $E_1, \ldots, E_n$,
 from (\ref{EQ:FHEV}) it holds  that
\[
P(E_i)=1, i=1,\ldots,n \; \Longleftrightarrow \; P(E_1 \cdots E_n)=1.
\] 
In the case where $P(E_i)=1, i=1,\ldots,n$, is coherent, 
 we can say that  the family $\F=\{E_1,\ldots, E_n\}$ is  \emph{p-consistent} and hence this property amounts to the coherence of the assessment $P(E_1 \cdots E_n)=1$. 
 This characterization still holds when replacing each $E_i$ by $E_i|H_i$; indeed, a family  $\F=\{E_1|H_1,\ldots, E_n|H_n\}$ is p-consistent when the assessment $P(E_i|H_i)=1, i=1,\ldots,n$, is coherent. Moreover, from (\ref{EQ:FHBOUNDS}) we obtain 
\[
P(E_i|H_i)=1, i=1,\ldots,n \; \Longleftrightarrow \; \prev[(E_1|H_1)\wedge  \cdots \wedge (E_n|H_n)]=1,
\] 
that is p-consistency of $\F$ amounts to the coherence of  the assessment $\prev[(E_1|H_1)\wedge \cdots \wedge(E_n|H_n)]=1$.
\\
(f) 
Given a p-consistent family $\F=\{E_1, \ldots, E_n\}$ and a further event $E_{n+1}$, the family $\F$ p-entails the event $E_{n+1}$ if and only if $P(E_i)=1, i=1,\ldots,n$, implies that $P(E_{n+1})=1$. The p-entailment of $E_{n+1}$ from $\F$ is equivalent to each one of the following properties (\cite{GiSa21A}): \\
$(i)$  $E_1 \cdots E_n \subseteq E_{n+1}$, that is  $E_1 \cdots E_n E_{n+1}=E_1 \cdots E_n$;\\
$(ii)$ the (indicator of the) conditional event $E_{n+1}|E_1 \cdots E_n$ is constant and coincides with 1.\\
The same characterization holds in the case of  conditional events. More precisely, given a p-consistent family  $\F=\{E_1|H_1, \ldots, E_n|H_n\}$ and a further conditional event $E_{n+1}|H_{n+1}$, we recall that  $\F$ p-entails $E_{n+1}|H_{n+1}$ if and only if $P(E_i|H_i)=1, i=1,\ldots,n$, implies $P(E_{n+1}|H_{n+1})=1$. This notion of p-entailment extends in the setting of coherence the notion introduced by Adams in \cite{adams75} by using \emph{proper} probability distributions. In (\cite[Theorem 18]{GiSa19}) it has been shown that the p-entailment of $E_{n+1}|H_{n+1}$ from  the  p-consistent family $\F$ can be characterized by  the following equivalent properties:\\
 $(i)$ $(E_1|H_1) \wedge \cdots \wedge (E_n|H_n) \leq E_{n+1}|H_{n+1}$; \\ 
 $(ii)$ $(E_1|H_1) \wedge \cdots \wedge (E_n|H_n) \wedge E_{n+1}|H_{n+1} = (E_1|H_1) \wedge \cdots \wedge (E_n|H_n)$. \\
Finally, 
by introducing a suitable notion of iterated conditional $(E_{n+1}|H_{n+1})\,\mid\,[(E_1|H_1) \wedge \cdots \wedge (E_n|H_n)]$, 
in \cite{GiSa21E} it has been shown that 
\begin{equation}\label{EQ:PEITER1}
	\F \; \mbox{p-entails} \; E_{n+1}|H_{n+1} \; \Longleftrightarrow \; (E_{n+1}|H_{n+1})\,\mid\,[(E_1|H_1) \wedge \cdots \wedge (E_n|H_n)] \, = \, 1 \,,
\end{equation}
that is  the p-consistent family $\F$ p-entails $E_{n+1}|H_{n+1}$ if and only if the iterated conditional $(E_{n+1}|H_{n+1})\,\mid\,[(E_1|H_1) \wedge \cdots \wedge (E_n|H_n)]$ is constant and coincides with 1. The characterization of the p-entailment for several inferece rules has been studied in \cite{GiPS20,GiSa21A}. 
\section{Conclusions}
In the first part of this work we illustrated the subjective approach to probability of de Finetti. We recalled the coherence principle for both the betting scheme and the penalty criterion, in the unconditional and conditional cases.
We showed the equivalence of the two criteria and    illustrated the geometrical interpretation of coherence. We also considered the notion of coherence  in the framework  of proper scoring rules.  We discussed the notion of conditional events in the trivalent logic of de Finetti, by focusing on the numerical representation of the truth-values. Then we considered some basic logical and probabilistic properties, valid for unconditional events, by checking their validity  in the setting of some  trivalent logics for some conjunction and disjunction pairs: $(\wedge_K, \vee_K)$,  $(\wedge_L, \vee_L)$, $(\wedge_B, \vee_B)$, $(\wedge_S, \vee_S)$. We verified that none  of these pairs satisfies all the properties; in particular, we observed that 
for these trivalent logics  the Fréchet-Hoeffding probability bounds are not satisfied. Then, we considered  our approach to compound conditionals in the setting of coherence and, based on the betting scheme, we have shown that conjunction and disjunction are suitable conditional random quantities.  Within this probabilistic logic, we verified that, when considering conditional events,  all the basic logical and probabilistic properties valid for unconditional events are preserved. In particular the Fréchet-Hoeffding probability bounds are  satisfied.  Moreover, we observed that the conjunction $E_1\cdots E_n$ of $n$ unconditional events coincides with the conjunction  $\C=E_1\wedge (E_2|E_1)\wedge \cdots \wedge (E_n|E_1\cdots E_{n-1})$. As a consequence, in our approach the prevision of $\C$ is the product  $P(E_1)P(E_2|E_1)\cdots P(E_n|E_1\cdots E_{n-1})$.
We also recalled the notions of p-consistency and p-entailment by illustrating their characterization by the notion of conjunction given in our approach.


\end{document}